\documentclass[11pt]{amsart}
\usepackage{amssymb}
\pdfoutput=1

 \title[CR Deformations]{The space of Cauchy-Riemann
structures on 3-D compact contact manifolds }

\date{June 30, 2010}

 \author{J. Bland  \and T. Duchamp}
	\address{John Bland\\ University of Toronto}
	\address{Tom Duchamp\\ University of Washington}

\thanks{The first author was partially supported by an NSERC
		  grant. The second author was partially supported by an
		  NSF grant.}



\numberwithin{equation}{subsection}

\newtheorem{lemma}[equation]{Lemma}
\newtheorem{proposition}[equation]{Proposition}
\newtheorem{theorem}[equation]{Theorem}
\newtheorem{corollary}[equation]{Corollary}

\theoremstyle{definition}

\newtheorem{definition}[equation]{Definition}

\newtheorem{remark}[equation]{Remark}

\renewcommand\Re{\mathrm{Re}}
\renewcommand\Im{\mathrm{Im}}

\newcommand\Id{\mathrm{Id}}

 \newcommand\G{\Gamma}
 \newcommand\Wk[1]{\Gamma^{#1}}

 \newcommand\Wkcont[1]{\Gamma_{cont}^{#1}(TM)}

 \newcommand\Wkcxcont[1]{\Gamma_{cont}^{#1}(T_{(1,0)}M)}

 \newcommand\WkCR[1]{\Gamma_{CR}^{#1}(TM)}

\newcommand\piH{\pi_{H}}

\newcommand\DeltaQ{\Delta_{R}}

\newcommand\diff[1]{\mathcal{D}^{#1}(M)}
\newcommand\cdiff[1]{\mathcal{D}_{cont}^{#1}(M)}

\newcommand\Thol{T_{(1,0)}M}

\newcommand\Hhol{H_{(1,0)}}
\newcommand\Hconj{H_{(0,1)}}

\newcommand\Hhat{\hat{H}}
\newcommand\wg{\wedge}

\newcommand\locdiff[1]{\mathcal{D}^{s}}

\newcommand\st{\; : \;}

\newcommand\norm[1]{\| #1\|}

\renewcommand\d{\partial}
\newcommand\pd[2]{\frac{\partial #1 }{ \partial #2}}
\newcommand\db{\bar{\partial}}
\newcommand\bdb{\db_{b}}

\renewcommand\a{\alpha}
\renewcommand\b{\beta}

\newcommand\g{\gamma}

\newcommand\ab{\bar{\alpha}}
\newcommand\bb{\bar{\beta}}
\newcommand\gb{\bar{\gamma}}

\newcommand\n{\eta}

\newcommand\w{\omega}
\newcommand\W{\Omega}

\newcommand\bw{\bar{\w}}

\newcommand\bW{\overline{W}}

\newcommand\str{{*}}

\newcommand\R{\mathbb{R}}
\newcommand\C{\mathbb{C}}

 \newcommand\Lie{\mathcal{L}}

 \newcommand\cE{\mathcal{E}}

 \newcommand\cH{\mathcal{H}}

\newcommand\cL{\mathcal{L}}

\newcommand\cN{\mathcal{N}}

\newcommand\cP{\mathcal{P}}

 \newcommand\cS{\mathcal{S}}
 
 \newcommand\cU{\mathcal{U}}
 \newcommand\cV{\mathcal{V}}

 \newcommand\inter{\mbox{\rule{2.0mm}{0.1mm}\kern0.0mm\rule{0.1mm}{3.0mm} }\,}

 \newcommand\lang{\ {{\vrule width .15in height .03pt
 depth .03pt}\!\!\!\!\!\joinrel{\raise3pt\hbox{$\backslash$}}}\ }

\newcommand{\XH}{X_{H}}
\newcommand\const[2]{C}

\newcommand{\FX}{F_{X}}            

\newcommand{\FPX}{F_{\Psi(X) }}    
\newcommand{\FPY}{F_{\Psi(Y) }}    

\newcommand{\Quad}{{Q}}   

\newcommand{\modeta}{\;\mathrm{mod }\; \n}

\newcommand{\Deform}{\,{\mathcal{D}ef}}

\newcommand{\apriori}{{\it a priori\/}}

\newcommand\bDb{{\bdb}}

\newcommand{\Hh}{H}
\newcommand{\Sh}{S}
\newcommand{\Qh}{Q}

\newcommand{\Ph}{P}

\newcommand{\Lh}{L}
\newcommand{\Nh}{N}

\newcommand{\cPh}{\cP}
\newcommand{\cQh}{\cH}
\newcommand{\cLh}{\cL}
\newcommand{\cNh}{\cN}
\newcommand{\cSh}{\cS}

\newcommand\rh{\rho}

\newcommand{\Harm}{\mathrm{H}^{1}}

\newcommand\range{\mathrm{range}}

\newcommand\HC{H_{\C}}           

\newcommand{\TMC}{T_{\C}M}

\newcommand\Forms[1]{\Omega^{(0,#1)}(M)}
\newcommand\FormsT[1]{\Omega^{(0,#1)}(\Thol)}
\newcommand\FormsH[1]{\Omega^{(0,#1)}(\Hhol)}

\newcommand{\formulae}{{formul\ae}}

\begin{document}
%
\begin{abstract}

We study the action of
the group of contact diffeomorphisms on CR deformations of
compact three-dimensional CR manifolds. Using anisotropic function
spaces and an anisotropic structure on the space of contact
diffeomorphisms, we establish the existence of local transverse slices 
to the action of the contact diffeomorphism group in the neighbourhood
of a fixed embeddable strongly pseudoconvex CR structure.

\end{abstract}

 \maketitle
 \section{Introduction}
\label{sec:intro}

Cauchy-Riemann manifolds arise naturally as the boundary of a bounded
domain $D \subset \C^{n+1}$. In this case, the Cauchy-Riemann structure is
simply that {\it residual} complex structure which is inherited from the
complex structure on $\C^{n}$. Local coordinates for $\d D$ are said to be
CR (for Cauchy-Riemann) if they are the restriction of holomorphic
coordinates in $\C^{n+1}$, and they define a conjugate CR tangent space for
$\d D$ in the same manner that the holomorphic coordinates on $\C^{n}$
define a conjugate holomorphic tangent space for $\C^{n+1}$. Intrinsically,
one can define the {\it Cauchy-Riemann} structure on $\d D$ by specifying
the space of conjugate CR tangent vectors in the same manner as one defines
the complex structure on $\C^{n+1}$ by specifying the conjugate holomorphic
tangent space. All questions which arise for abstract complex structures on
a smooth manifold are equally valid for Cauchy-Riemann manifolds: for
example, the embeddability and local embeddability (or the existence of
holomorphic (CR) coordinates, or how many structures exist up to
equivalence.

The significance of generalizing from complex structures on manifolds to
studying Cauchy-Riemann structures can easily be seen from the following
considerations. When $D$ is a bounded domain in $\C^{n+1}, n \ge 1$, then
holomorphic functions on $D$ which extend smoothly to $\d D$ restrict to
$\d D$ as CR functions; on the other hand, a slight generalization of
Hartog's phenomenon in several complex variables states that CR functions
on $\d D$ extend uniquely to $D$ as holomorphic functions; that is, $\d D$
with its Cauchy-Riemann structure completely determines $D$ with its
complex structure.  On the other hand, if we generalize to $\Sigma$ a
complex analytic space with an isolated singularity at $p \in \Sigma$, then
the boundary of a small neighbourhood of $\Sigma$ inherits a smooth Cauchy
Riemann structure whereas the space $\Sigma$ is singular.  On the basis of
this observation, Kuranishi proposed \cite{Ku} to study the deformation
space for isolated singularities by studying the deformation space for
Cauchy-Riemann structures on the boundary of the neighbourhood, a smooth
compact manifold.

A case of particular interest is that in which the domain $D$ is strongly
convex (more generally, strongly pseudoconvex). In this case, the boundary
admits a natural family of positive definite metrics which are adapted to
the CR structure, and play much the same role that K\"ahler metrics play in
complex geometry. One consequence of particular importance is that when $M$
is compact, strongly pseudoconvex and $n \ge 2$ (so $\dim M \ge 5$), then
$M$ is embeddable.  This is definitely not the case when $n=1$, and this
case has many deep and interesting features which have yet to be fully
understood.

In this paper, we fix a smooth compact underlying manifold, and study the
space of CR structures on the manifold up to equivalence. In particular, we
study the local deformation theory for the space of CR structures, and the
local action of the contact diffeomorphism group on the space of such
structures. Although for much of the paper we set up the machinery to work
in arbitrary dimensions, our main interest is in the three dimensional
case, and we restrict our attention to this case in the latter sections of
this paper This was largely a matter of expedience, since in higher
dimensions integrability factors play a role, and require the introduction
of new operators and significantly different treatment than in the three
dimensional case.

Most of the results in this paper rely heavily on \cite{BD1} in which
we developed the machinery to do analysis on contact manifolds using
intrinsically defined anisotropic functions spaces.

The outline of the paper is as follows. In Section 2, we give a quick
review of strongly pseudoconvex Cauchy-Riemann structures and the relevant
deformation theory.  In Section 3, we define the weighted or anisotropic
function spaces in which we will work, and recall the results from
\cite{BD1} on the space of weighted contact diffeomorphisms which we
will need throughout the remainder of the paper.  The inclusion of these
two sections is to fix notation and to help make the paper
self-contained. In Section 4, we study the action of contact
diffeomorphisms on CR structures, computing both the linear and the fully
nonlinear action; it is also in this section that we introduce the notion
of complex contact vector fields, and explain their relation to the
symmetry group. In Section 5, we collect results on homotopy operators for
the $\bdb $-complex on compact CR manifolds and adapt them to our
particular situation; we also indicate how to split complex contact vector
fields into real contact vector fields and a transverse vector
field. Section 6 contains the main results of the paper. In this section,
we obtain normal forms for CR structures under the action of the group of
contact diffeomorphisms with sharp regularity results. This is accomplished
in two steps: first we obtain a weak normal form with a loss of regularity,
and then using \apriori\ estimates we recover the lost regularity. It is
believed that this approach to studying the action of infinite dimensional
symmetry groups on underlying structures is new, and may have applications
in other situations.

Earlier results in this direction were obtained in \cite{CL} and and
\cite{B}.  The main idea in both papers was to study the linearized action,
and to construct appropriate function spaces in which one can solve the
linearized equation with good estimates.  Since the $\bdb$--operator
appears in the linearized equation, the anisotropic function spaces appear
naturally. In \cite{CL}, they avoided using the anisotropic spaces by
working in the Nash Moser category; they obtained a transverse slice for
smooth CR~structures. In \cite{B}, we restricted our attention to the case
of the standard $S^3 \subset \C^2$, and used explicit information to
construct an anisotropic Hilbert space structure on contact diffeomorphisms
near the identity; the description of transverse slices follows easily from
the linearized analysis. However in \cite{B}, the action described for the
contact diffeomorphism group was incorrectly asserted to be $C^{1}$, a
necessary condition to apply the inverse function theorem in Banach spaces
and obtain the transverse slices; a modified action is used in Section~
\ref{sec:slice} of the current paper to correct this error. With this
modification and the generalization of the weighted function space
structure for contact diffeomorphisms to arbitrary compact contact
manifolds (see \cite{BD1}), we are now able to obtain local transverse
slices to the action of the contact diffeomorphism group on the space of CR
structures for an arbitrary compact embeddable strongly pseudoconvex three
dimensional CR structure.

\subsection{Notation}
\label{notation-2}

Throughout the paper, $M$ will denote a smooth compact $2n+1$
dimensional manifold equipped with a fixed contact distribution $H \subset TM$
and a fixed contact one form $\n$.
As usual,   $TM$ and $T^*M$ denote the tangent and
cotangent bundles of $M$, respectively, $\Lambda^p M$ denotes
the $p$-th exterior power of $T^*M$,  $\W^{p}(M)$ the space of
smooth $p$-forms on $M$,  $\Lie_{X}\beta$  the Lie derivative of the
form $\beta$ with respect to the vector field $X$, and $X\inter \beta$
interior evaluation.

We give $M$ a fixed Riemannian metric $g$ compatible with $\eta$ (see
Equation~\eqref{eq:g} for details), and let $|X|$ denote the
norm of the tangent vector $X$ with respect to $g$, and we let $\exp:TM\to
M$ denote the exponential map of the $g$.

We let
\[
   \piH : T^{*}M \to H^{*} 
\]
denote the projection map.  The {\em characteristic}
(or {\em Reeb}) vector field $T$ is the unique vector field satisfying the
conditions $T\inter \n = 1$ and $T\inter d\n = 0$. We can then identify the
dual contact distribution with the annihilator of $T$, i.e.
\[
     H^{*}
= \{ \b \in T^{*}M \st   T \inter \b = 0 \} \subset T^{*}M \,;
\]
more generally
\[
          \Lambda^{p}H^{*} = \{ \b\in \Lambda^{p}(M) \st T \inter \b = 0 \}\,,
\]
and we have the identity
\begin{equation}
\label{eq:piH-formula}
           \piH(\b) = T \inter ( \n \wedge \b) \,.
\end{equation}

We endow $\R^{2n+1}$ with
the contact structure defined by the one-form
\[
	 \n_{0} =  dx^{2n+1} -  \sum_{j=1}^{n} x^{n+j} dx^j \,,
\]
where  $(x^1,\dots,x^n,x^{n+1},\dots,x^{2n},x^{2n+1})$ are the standard
coordinates on $\R^{2n+1}$, and we let $dV_0$ denote the standard volume
form:
\[
        dV_0 = \frac{1}{n!} \n_0\wedge (d\n_0)^n\,.
\]
We  denote the contact distribution of $\n_{0}$ by  $H_0 \subset T\R^{2n+1}$
and we set
\[
  T_{0} = \pd{}{x^{2n+1}}\,, X_j = \pd{}{x^j} + x^{n+j}\, \pd{}{x^{2n+1}}\,, \mathit{ and } X_{n+j} = \pd{}{x^{n+j}}
\,, 1 \leq j \leq n \,.
\]
Observe that the collection $\{X_j, 1\leq j\leq 2n
\}$ is a global framing for $H_0$.
Note also that the 1-forms
\[
	 \n_0,  dx^{j},\, dx^{n+j},\, 1 \leq j \leq n,
\]
are the dual coframe to $T_0$, $X_j$, $X_{n+j}$, $1 \leq j \leq n$.

Let $f = (f^1,\dots,f^m)$ be a smooth, $\R^m$-valued function defined on the
closure of a domain $D \Subset\R^{2n+1}$. We define
 \[
	X_{I} f = \begin{cases}
		     X_{i_1}X_{i_2}\dots X_{i_t}f &\text{for $t>0$}\\
		     f &\text{for $t=0$} \,,
		  \end{cases}
\]
where we have introduced the multi-index notation $ I=(i_1,\dots,i_t)$,
$1\leq i_j\leq 2n$ and $X_{I}f = (X_{I}f^1,\dots,X_{I}f^m)$. (For $t=0$,
$I$ denotes the \emph{empty index} $I=()$.)  The integer $t$ is called the
\emph{order} of $I$ and written $|I|$.

\begin{remark}
\label{rem:local-coords}
We will often have to work in local coordinates adapted to the contact
structure on $M$.
An \emph{adapted coordinate chart} for $M$ is a chart $\phi:U \to
\R^{2n+1}$ for which $\n = \phi^*\n_{0}$. It follows that $\phi_* T = T_0$
and $\phi_* H = H_0$. 
An \emph{adapted atlas} consists of the following data:
a fixed
 finite open cover $V_{\ell} \Subset U_{\ell}$, $\ell=1,2,\dots,m$
 and an atlas  $\{\phi_{\ell}: U_{\ell} \to \R^{2n+1}\}$, consisting of adapted
coordinate charts. We set $D_{\ell} = \phi_{\ell}(V_{\ell})$.
By compactness of $M$ and Darbourx's Theorem for contact structures
\cite[page 362]{Arn:1978}, $M$ has an adapted atlas.  
We shall fix once and for all an adapted atlas and a partition of unity $\rho_{\ell}$
 a partition of unity subordinate to   $\{ V_{\ell} \}$.

\end{remark}

If $\mathcal{F}:\mathcal{A} \to \mathcal{B}$ is a map between Banach
 spaces, with norms $\|\,\cdot\,\|_{\mathcal{A}}$ and
  $\|\,\cdot\,\|_{\mathcal{B}}$, respectively, then the expression
 \[
		     \|\mathcal{F}(f)\|_{\mathcal{B}} \prec
		     \| f \|_{\mathcal{A}}
 \]
 means that there is a constant $C>0$ such that
 $\|\mathcal{F}(f)\|_{\mathcal{B}} \le C \| f \|_{\mathcal{A}}$ for all $f \in
 \mathcal{A}$.


\section{CR~structures}
\label{sec:CR-structures}

\subsection{Deformation theory of CR~structures}

We begin with a quick review of the deformation theory of CR~structures as
presented in the paper of Akahori, Garfield, and Lee
\cite{Aka-et-al:2002}. See also \cite[Section~16]{BD2}, where the special
case of the deformation theory of $S^{2n+1}$ is studied using a similar
framework.

\begin{definition}
\label{def:CR-structure}
Let $M$ be a $2n+1$ dimensional manifold.
A (rank $n$) \emph{Cauchy-Riemann structure (CR-structure)} on
 $M$ is a rank $n$ complex subbundle $H_{(1,0)} \subset
T_{\C}M$ of the complexified tangent bundle of $M$ such
that
\begin{enumerate}
\item[(i)]
$H_{(1,0)} \cap \overline {H_{(1,0)}} = \{0\}$,
\item[(ii)]
the integrability condition is satisfied:
\[       [ \Gamma^{\infty}(\Hhol) , \Gamma^{\infty}(\Hhol) ] \subset
      \Gamma^{\infty}( \Hhol) \, .
\]
\end{enumerate}
The bundle $\Hhol$ is called the \emph{holomorphic tangent bundle} of
the CR-structure.
As usual, we let $\Hconj$ denote the conjugate bundle
$\overline{H_{(1,0)}}$. The transversality condition (i) implies that
$H_{\C} =H_{(1,0)} \bigoplus \Hconj \subset T_{\C}M$ has complex
codimension one.
 \end{definition}

\begin{remark} 
We recall that when $n=1$, the bundle $\Hhol$ is a complex line bundle,
and condition (ii) is automatic. To see this let $Z$ be section of
$\Hhol$ that does not vanish on an open set $U \subset M$. Then
since $[Z,Z] = 0$,
\[
     [f Z, g Z] = \left( f Z(g) - g Z(f) \right) Z  
\]
for any two sections, say $X = f Z $ and $Y = g Z$ of $\Hhol$.
\end{remark}

Two CR-structures $\Hhol$ and $\Hhat_{(1,0)}$ are said to be
\emph{equivalent} if there is a diffeomorphism $F:M\to M$ such that
$F_{*} \Hhol = \Hhat_{(1,0)}$.  We are only interested in CR-structures
up to equivalence.  

Observe that $H_{\C}$ is the complexification of a real
codimension one subbundle $H \subset TM$ consisting of vectors of
the form $X+\overline{X}$, $X\in \Hhol$. Let $\n$ be a real one-form dual
to $H$.
The CR-structure $\Hhol$ is said to be \emph{strongly pseudoconvex} if
$-i\, d\n(X, \overline{X}) > 0$ for all non-zero $X \in \Hhol$. In this
case, $\n\wg(d\n)^n$ is a nowhere vanishing $(2n+1)$-form.  In other
words,
$(M,H)$ is a contact manifold and $\n$ is a contact one-form.

The most common examples of CR-structures are those arising from
domains in $\C^{n+1}$.  Let  $D = \{z\in\C^{n+1} \st \rh(z) <1 \}$ 
be a smoothly bounded domain in $\C^{n+1}$ with connected boundary, where 
$\rh$ is a smooth nonnegative function defined on a neighbourhood of
$\overline D$, and $d \rh \ne 0$ on $\d D$.  The boundary $\d D$ is a
CR-manifold for which the
holomorphic tangent bundle is the intersection of the complexified
tangent bundle of $\d D$ with the holomorphic tangent bundle of
$\C^{n+1}$, and if the pullback to $\d D$ of the one-form $i\db \rh$
is
a contact form, then it is strongly pseudoconvex\footnote{The fact
that $D$ is bounded forces the Levi form to be positive at some point 
on $\d D$, hence by the non-degeneracy of $d \n$ everywhere on the connected
manifold $\d D$.}.

We will assume, henceforth, that $M$ is a contact manifold
 equipped with a fixed strongly pseudoconvex CR-structure  
$\Hhol \subset T_{\C}M$ such that $\HC$
is the complexification of the contact distribution of $M$. 
We shall refer to this CR-structure as the \emph{reference CR-structure} on
$M$.

The reference CR~structure determines an endomorphism $J : H \rightarrow H$
satisfying the condition $J^2 =
-Id$, which in turn defines a Riemannian metric $g$ by the formula
\begin{equation}
\label{eq:g}
     g( X, Y) = \n (X) \n(Y) + d\n (X, J Y)\,.
\end{equation}
The metric $g$ is said to be \emph{adapted to the CR~structure}. We let
$\exp: TM \rightarrow M$ denote the  exponential map of $g$.
Objects associated to any other CR-structure on $M$ will be decorated with hats.
Two strongly pseudoconvex CR-structures on $M$ are said to be
\emph{isotopic} if they can be connected by a smooth 1-parameter family
of strongly pseudoconvex CR-structures. We consider
only strongly pseudoconvex CR-structures which are isotopic to the
reference CR-structure.

\subsection{Representation by Deformation Tensors}
Every CR-structure that is isotopic to the reference one can be
represented by a \emph{deformation tensor} that takes values in
$\Hhol$.
The proof of this fact relies on a theorem of John Gray \cite{gray:1959a}
which states that isotopic contact structures on a compact manifold  
are equivalent.

\begin{theorem}[Gray]
Let $\n_t$ be a differentiable family of contact forms on a compact
$2n+1$ dimensional manifold $M$. Then there is a differentiable family of
diffeomorphisms $F_t : M \to M$ and a family of non-vanishing functions
$p_t$ such that
\[
	    F_t{}^{*} (\n_t)  =  p_t \n_0\,.
\]
\end{theorem}

\begin{corollary}
\label{cor:deform-tensor}
Every strongly pseudoconvex CR-structure on $M$ that is isotopic to
the reference one is CR-equivalent to one of the form $\Hhat_{(1,0)}$ 
where
\begin{equation}
\label{eqn:Hhat}
     \Hhat_{(0,1)} = \{   X  - {\phi}(X) \st X \in \Hconj \}
\end{equation}
and $\phi: \Hconj \to \Hhol$ is a map of complex vector bundles, called
the \emph{deformation tensor} for $\Hhat_{(1,0)}$.
\end{corollary}

\begin{proof} 
The fact that the CR-structure is equivalent to one satisfying the
inclusion
relation $\Hhat_{(0,1)} \subset \HC$ follows immediately from Gray's
theorem.  
Thus, there is a family $\Hhat_{(0,1)}(t)$, $t \in [0,1]$, joining
$\Hconj$ to
$\Hhat_{(0,1)}$.  For $t$ small, it is clear that there are bundle maps
$\phi(t)$ such that $\Hhat_{(0,1)}(t)$ is the graph of $-{\phi}(t)$.
The integrability conditions for CR-structures 
imply that $\phi(t)$ satisfies certain symmetry properties, and when 
combined with the transversality condition they imply an \apriori\ bound
on the size of $\phi(t)$, from which the result follows.  

We explain in brief. Choose a local basis 
$Z_{\a}$ for $\Hhol$, and let $i \n [Z_{\a}, Z_{\bb}] =
-i d\n (Z_{\a}, Z_{\bb}) =  h_{\a \bb}$ 
define the Levi form.\footnote{
Here, and for the remainder of this
  section, we  employ the Einstein summation conventions, with Greek
  indices ranging from $1$ to $n$, and the conventions for raising and 
lowering indices by contraction with the hermitian form $h_{\a
\bar\delta}$ and its inverse, with 
$\phi_{\bb}^{\g} h_{\g \bar\delta} = \phi_{\bb \bar\delta}$ .}
Then integrability implies in particular that
$i \n [Z_{\bb}-\phi^{\a}_{\bb} Z_{\a},
Z_{\bar\delta}-\phi^{\g}_{\bar\delta} 
Z_{\g}] =0$. Since 
$\n [Z_{\a}, Z_{\g}] = \n [ Z_{\bb}, Z_{\bar\delta}] =0$, it
follows  that
$i \n [-\phi^{\a}_{\bb} Z_{\a}, Z_{\bar\delta}] + i \n
[Z_{\bb}, -\phi^{\g}_{\bar\delta} Z_{\g}]
= \phi^{\g}_{\bar\delta} h_{\g \bb} -\phi_{\bb}^{\a} h_{\a \bar\delta}=0$;
this is the symmetry condition 
$\phi^{\g}_{\bar\delta} h_{\g \bb} =\phi_{\bb}^{\a} h_{\a \bar\delta}$. It
follows that the operator $\phi \circ \bar\phi : \Hhol
\rightarrow \Hhol $  has non-negative
eigenvalues since 
\[
 i \n [(\phi \circ \bar\phi) Z_{\a}, Z_{\bar \mu}] = 
(\phi \circ \bar\phi)_{\a}^{\delta} h_{\delta \bar \mu} = 
(\bar\phi)_{\a}^{\bar \g} \phi_{\bar \g}^{\delta}  h_{\delta \bar \mu}
 = h^{\bar \g \b} (\bar\phi)_{\b \a} \phi_{\bar \g \bar \mu} 
\]
is hermitian positive semi-definite.

Next note that the transversality condition for CR structures 
(Definition~\ref{def:CR-structure}(ii)) 
implies that none of the
eigenvalues of $\phi \circ \bar \phi$ can be equal to one. Indeed, 
suppose to the  contrary that 
$(\bar\phi)_{\a}^{\bar \g} \phi_{\bar \g}^{\delta} - \delta_{\a}^{\delta}$
is a degenerate matrix. Then there exists $v^{\a}$ such that 
$v^{\a} (\bar\phi)_{\a}^{\bar \g} \phi_{\bar \g}^{\delta} =
v^{\delta}$, from which one obtains the relation
\[
   v^{\a}(Z_{\a}-(\bar \phi)^{\bb}_{\a} Z_{\bb}) =  
   v^{\a}(\bar \phi)_{\a}^{\bar \b} \phi_{\bar \b}^{\delta} Z_{\delta}
-  v^{\a}(\bar \phi)^{\bb}_{\a} Z_{\bb} =
-  v^{\a}(\bar \phi)^{\bb}_{\a} \left( Z_{\bb} - \phi_{\bar \b}^{\delta}
Z_{\delta}\right);
\]
that is, the transversality condition is violated for the
subspace $\Hhat_{(1,0)}(t)$ and its conjugate.

Since $\phi \circ \bar \phi$ is isotopic to the zero map by assumption,
has
non-negative eigenvalues, and $(\phi \circ \bar \phi - I)$ is
nondegenerate, it follows that the eigenvalues of the
operator $(\phi \circ \bar\phi)$ are bounded between $0$ and $1$,
which  implies the norm condition.
(See \cite[page
83]{BD2} where a similar argument is given.)
\end{proof}

\begin{remark}
The choice to refer to the map $\phi: \Hconj \to \Hhol$
as the deformation tensor (rather than the conjugate map) is
consistent with the deformation theory for complex structures, and has 
the advantage that $\phi$ may be thought of as a ``vector-valued
$(0,1)-$form'', thus fitting naturally within a $\db -$complex
(or in this case, a $\bdb -$complex).
\end{remark}

In light of Corollary~\ref{cor:deform-tensor}, we  identify the space of
CR~structures with the subset of the space
of $\Hhol$-valued $(0,1)$-forms.
Specifically, if $\w^{\a}$ is a local coframe of $H^{(1,0)}$ with dual  
frame $Z_{\a}$ of $\Hhol$  such that
\[
	      d\n = \frac{i}{2} \delta_{\a\bb} \w^{\a}\wedge \w^{\bb}\,,
\] 
then the \emph{CR~deformation tensor} can be written as
\[
      \phi = \phi_{\bb}^{\a} \, \w^{\bb}\otimes Z_{\a} \, ;
\] 
it uniquely determines the space of $(0,1)$-vectors for its
corresponding CR~structure as   
the subspace of $\HC$  annihilated by
 the one-forms 
\begin{displaymath}
	 \hat{\w}^{\a} := \w^{\a} + \phi^{\a}_{\bb} \w^{\bb}  \; .
\end{displaymath}
The space of all smooth deformation tensors is given  by
\begin{equation}
\label{eq:def-deformation}
    \Deform =  \FormsH{1} \simeq \Wk{\infty}\left( H_{(0,1)}\otimes\Hhol
          \right) \,.
\end{equation}

\subsection{The deformation complex}
Each deformation of a CR~structure can be expressed as a $\Hhol$-valued
(0,1)-form. In \cite{Aka:1981}, Akahori studied CR~deformations by
developing the Hodge theory of a certain complex of vector-valued
forms. A
similar complex was studied in \cite{BD2} and used to show that
CR~deformations of the standard CR~structure on $S^{2n+1}$ can be
parameterized by \emph{complex Hamiltonian vector fields}.

The space of smooth \emph{forms of type $(0,q)$}, written
$\Forms{q}$, is the space of sections of the bundle $\Lambda^q
H^{(0,1)}$,
where $H^{(0,1)}$ denotes the dual bundle of the complex vector bundle
$\Hconj$. By the integrability condition for the CR structure, the exterior
differential operator $d$ naturally induces an operator
\[
     \bdb: \Forms{q} \to \Forms{q+1} \, .
\]

Set $\Thol = \TMC / \Hconj$, where
$\TMC = TM \otimes_{\R} \C$ is the complexified tangent bundle of $M$,
and
let $\pi_{(1,0)}:\TMC \to \Thol$ denote the quotient map. The space of
\emph{$\Thol$-valued forms of type $(0,q)$}  is the space of
homomorphisms
of complex vector bundles
\[
	 \FormsT{q} =  \G^{\infty}
    \left( \mathrm{Hom}_{\C}( \Lambda^q \Hconj M,\Thol)\right) \,.
\]  

By virtue of the integrability condition
(Definition~\ref{def:CR-structure}(ii)), the operator $\bdb$ extends to an
operator on the space of $\Thol$-valued forms \cite{BlEp,BuMi},
which by abuse of notation we again denote by $\bDb$:
\begin{equation}
\label{eqn:vect-dbar}
     \bDb: \FormsT{q}   \to \FormsT{q+1} \,.
\end{equation}
This operator is characterized by the following  properties:
\begin{subequations}
\begin{align} \bDb^2 &= 0 \,;\\
\label{eq:bDbX}
 \bDb(X)(\overline{Z}) &= \pi_{(1,0)} [\overline{Z}, X ]\,,
\intertext{for $X \in \FormsT{0} = \Gamma(\Thol)$ and $\overline{Z} \in
\Gamma(\Hconj M)$;}
    \bDb( \a \wedge \b) &= (\bdb\a) \wedge \b + (-1)^{q_1} \a \wedge
    \bDb\b
    \,,
\end{align}
for $\a\in\Forms{q_1}$ and  $\b\in\FormsT{q_2}$.
\end{subequations}
\begin{remark}
\label{rem:lifted-operator}
The operator defined by equation~(\ref{eq:bDbX}) further lifts to an
operator 
\[
     \bDb: \G(T_{\C}M )  \to \FormsT{1} \,
\]
via the formula
\[
           \bDb X := \bDb (\pi_{(1,0)}X) \,.
\]
In particular, $\bDb X$ is well-defined in the special case where $X$ is a
real vector field.
By abuse of notation we again denote the lifted
operator by
$\bDb$.  
\end{remark}

\begin{remark}
\label{rem:embedded}
When $(M, \n, \Hhol)$ is embedded, it bounds a strongly
pseudoconvex complex space $\Sigma$, and there
is a natural identification between $\Thol$ and the restriction of the
holomorphic tangent bundle from $\Sigma$. In this case, $\bDb$
is naturally identified with the restriction of the $\db$
operator to the boundary.
\end{remark}

\section{The group  of Folland-Stein contact diffeomorphisms}
\label{sec:top-group-str}

In the previous section, we showed that every CR~structure isotopic to a
reference CR~structure can be represented by a deformation tensor. In
Section~\ref{sec:action}, we study the action the group of contact
diffeomorphisms on to space of CR~structures. In this section, we recall
the results from \cite{BD1} that we need.  Details can be found
in \cite{BD1}.

\subsection{Folland-Stein spaces}

We begin by recalling the anisotropic function spaces $\Wk{s}(M)$ on
$M$, introduced by Folland and Stein in \cite{FS}, and their
generalizations. These spaces are the natural ones in which to work in
order to obtain sharp estimates for the various operators which will
arise. 


 Consider an open domain $D \Subset \R^{2n+1}$.
The {\em Folland-Stein space} $\Wk{s} = \Wk{s}(D)$ is the Hilbert space
completion of the set of smooth functions on $\overline{D}$ with respect to
the
inner product
\[
   (f,g)_{D,s} :=  \sum_{0 \leq |I| \leq s} \int_{D}
      |X_I f|\, |X_I g| dV \, ,
\]
with associated norm written  $\|f\|_{D,s} = \sqrt{(f,f)}$, where
 $X_I$ and $dV$ are as in Section~\ref{notation-2}. 
Let $\Wk{s}(D,\R^m)$ denote the closure of the smooth $\R^m$
valued functions on $\overline{D}$ with inner product
\[
(f,g)_{D,s} = \sum_{j=1}^{m} (f^j,g^j)_{s}
\]
for smooth functions  $f = (f^1,\dots,f^m)$ and $g = (g^1,\dots,g^m)$.


Let $(M, \n)$ be a smooth compact contact manifold, and let
$\{\phi_{\ell}:\overline{U}_{\ell} \to \R^{2n+1}\}$ be an adapted atlas as in
Section~\ref{notation-2}). 
A function $f:M \to\R$ is said to be a {\em $\Wk{s}$ function} if the
functions $f_{\ell} = f \circ \phi_{\ell}^{-1}$ lie in $\Wk{s}(D_{\ell})$ for
all
$\ell$. For functions $f, g \in \Wk{s}(M)$, we 
define the inner product
\[
   (f,g)_{s} := \sum_{\ell} (\rho_{\ell} f_{\ell},\rho_{\ell}
g_{\ell})_{D_{\ell},s} \, .
\]
Similar definitions hold for $\Wk{s}(M,\R^m)$, $f,g \in
\Wk{s}(M,\R^m)$.
The definition of the function spaces is independent of the choice of 
adapted atlas and the local framings $X_I$ and $dV$. Although the 
definition of the inner products depend upon the choices involved, 
different choices lead to equivalent norms.

Let $F:M \to \widetilde{M}$ be a $C^1$-map from $M$ into a smooth
$\widetilde{m}$-dimensional manifold $\widetilde{M}$.  Choose an adapted
atlas $\{ (\phi_{\ell},U_{\ell},V_{\ell})\}$ for $M$ and a smooth atlas
$\{\widetilde{\phi}_{\ell}:\widetilde{U}_{\ell}\to \R^{\widetilde{m}}\}$
for $\widetilde{M}$ such that
\[
F(U_{\ell}) \subset \widetilde{U}_{\ell},
\text{ and }
F_{\ell}(\bar{D}_{\ell}) \subset \widetilde{D}_{\ell}
\]
for all $\ell$,
where
$\displaystyle
 F_{\ell} = \widetilde{\phi}_{\ell} \circ F \circ \phi_{\ell}^{-1}: 
\phi_{\ell}\left(U_{\ell}\right) \to  \R^{\widetilde{m}}
$
and $\widetilde{D}_{\ell}\Subset \phi_{\ell}(\widetilde{U}_{\ell})$ is a
collection of open domains such that
$\{\widetilde{\phi}_{\ell}^{-1}(\widetilde{D}_{\ell})\}$ covers
$\widetilde{M}$.  The map $F$ is said to be a {\em $\Wk{s}$ map} if
$F_{\ell}$ restricts to an element $F_{\ell} \in
\Wk{s}(D_{\ell},\R^{\widetilde{m}})$ for all $\ell$.  It is not difficult to
show that the notion of $\Wk{s}$ map is independent of the choice of
atlases and that $F_{\ell}$ restricts to an element in $\Wk{s}(D)$ for any
open set $D\subset\subset \phi_{\ell}\left(U_{\ell}\right)$.

Let $\Wk{s}(M,\widetilde{M})$ for $s \geq n + 4$ denote the topological
space of $\Wk{s}$ maps between $M$ and $\widetilde{M}$. 
The restriction $s \geq n+4$ ensures that the maps are
$C^1$. More generally, consider a smooth
fibre bundle $\pi: P \to M$, with base a compact contact manifold. The
space $\Wk{s}(P)$ of $\Wk{s}$ sections of $\pi$ is defined in the
obvious way by choosing an adapted atlas for $M$ such that
$\pi^{-1}(U_{\ell}) \to U_{\ell}$ is trivial for all $\ell$ and requiring
the local coordinate representations of sections to be $\Wk{s}$ maps
from $U_{\ell}$ into the fiber of $\pi$. (See \cite{BD1} for details.)


\subsection{The smooth manifold of Folland-Stein diffeomorphisms}
\label{sec:param}
\label{subsec:exp}

Let $\diff{s} \subset \Wk{s}(M,M)$ denote the space of $\Gamma^s$
diffeomorphisms of $M$.  We showed in \cite{BD1} that $\diff{s}$ is an open subset of
$\Wk{s}(M,M)$ for all $s \geq 2n+4$.
Let $\cdiff{s} \subset \diff{s}$ denote the subspace of $\Wk{s}$ contact
diffeomorphisms of $M$.  In \cite{BD1}, we obtained a local coordinate
chart for contact diffeomorphisms in a neighbourhood of the identity, and we
showed that $\cdiff{s}$ is a topological group with respect to composition,
provided that $s \geq 2n+4$.  

More precisely, let $g$ be a metric  adapted to the
contact structure such as the one constructed
in the Section~\ref{notation-2}.
The exponential map induces various maps between $\Wk{s}$ spaces that we
need to parameterize contact diffeomorphisms. If $X$ is a vector field,
we
use the notation $\FX$ to denote the map
\begin{equation}
\label{eq:defFX}
  \FX := \exp \circ X : M \to M \,.
\end{equation}
Recall that because $M$ is compact, the map $\FX$ is a diffeomorphism for
$X$ sufficiently small.  The following proposition summarizes various
smoothness properties of the maps that we need to construct our local
coordinate charts for contact diffeomorphisms.

\begin{proposition}
\label{prop:exp}
Let $\Wk{s}(TM)$ denote the space of $\Wk{s}$ sections of $TM$.
For $s \geq (2n+4)$, the map
\[  \exp : \Wk{s}(TM) \to \Wk{s}(M,M) \,:\,
		    X \mapsto \FX = \exp \circ X 
\]
is smooth.  
Moreover, there is a neighbourhood $\cU \subset \Wk{2n+4}(TM)$
such that $\FX$ is in $\diff{s}$ for all $X \in \cU^s$ and all $s\geq 2n+4$, where
$\cU^s := \cU\cap \Wk{s}(TM)$;
and the restriction
\[
   \exp : \cU^s \to \diff{s}
\]
is a homeomorphism from $\cU^s$ to a neighbourhood of the identity
diffeomorphism.

\end{proposition}

In general, the diffeomorphism
$\FX$ of Proposition~\ref{prop:exp}
will not be a contact diffeomorphism. However in \cite{BD1}, we
showed that the subset of $\cU^s$ for which it is a contact diffeomorphism is
smoothly parameterized by the set of contact vector fields in a
neighbourhood of the zero section. As shown in \cite{BD1}, this implies that the 
space of $\Gamma^s$ contact diffeomorphisms is a smooth Hilbert manifold.

We now introduce some notation that will be necessary to express the
sharp estimates used later in the paper. 
Choose an adapted atlas $\phi_{\ell}: U_{\ell}\to \R^{2n+1}$
for $M$ and a collection of open sets $V_{\ell}\Subset U_{\ell}$ covering $M$
as in Section~\ref{notation-2} and let $\rho_{\ell}$ be a partition of unity
subordinate to $\{V_{\ell}\}$.  By compactness of $M$, there is a
constant $c>0$ such that $\exp(x,X) \in U_{\ell}$ for all $x \in
\overline{V_{\ell}}$, all $X\in TM_x$, with $|X|<c$, and all $\ell$.
Let $X$ be a $C^1$ vector field with $|X|<c$. 

Fix a chart, say $\phi_{\ell}$, and
set $U = U_{\ell}$ and $V =V_{\ell}$. To simplify notation, we adopt the
Einstein summation conventions, letting Roman indices range from $1$ to
$2n+1$. As explained in \cite{BD1}, by the second order Taylor's formula
with integral remainder, there exist smooth functions $B^{k}_{ij}(x,X)$ (locally
defined) on $TM$ such that
\begin{equation}  
\label{exp-map}  \FX^{k} := \exp^{k}(x,X) = x^{k} +  X^{k} + B^{k}_{ij}(x,X)X^{i}X^{j}\, . 
\end{equation} 
A standard computation using Equation~\eqref{exp-map} then yields the
following expansion for the pull-back of a $q$-form by $\FX$.

\begin{lemma}[\cite{BD1}]
\label{lem:exp-q-form}
Let $\psi$ be a smooth $q$-form on $M$ and
choose a coordinate patch $U = U_{\ell}$, with $V = V_{\ell} \Subset U$. Let
$c>0$ be chosen so that $\exp(x,X) \in U$ for all $x \in \overline{V}$ and
all $X \in T_xM$ with $|X|<c$.
Then  there are (locally defined) smooth fibre bundle maps
\[
    \Quad^1_{ij}: \left.   BM\right|_{V} \to \left.\Lambda^{q}M \right|_{V}
\text{ and }
    \Quad^2_{ij}: \left. BM\right|_{V} \to \left.\Lambda^{q-1}M \right|_{V} \,,
\]
where  $BM  = \{ X \in TM \st |X| < c\}$,
such that  for any $C^1$ vector field $X:M \to BM \subset TM$ the equation
\[
	\FX^{*} \psi = \psi + \Lie_{X} \psi
       +  \Quad^1_{ij}( X )\,  X^{i}X^{j}
       +  \Quad^2_{ij}( X ) \, \wedge X^{i} d X^{j}
\]
is satisfied  on all of $V$.
\end{lemma}
Henceforth, we will use the notation
\begin{equation}
\label{eqn:calR-notation}
    \Quad_{\psi}(X) := \FX^*(\psi) -\psi - \Lie_X\psi  
\end{equation}
to denote the non-linear part of the pull-back $\FX^*\psi$.
 The lemma states that in local coordinates
\begin{equation}
\label{eq:Quad-local}
   \Quad_{\psi}(X) =   \Quad^1_{ij}(X) \, X^{i}X^{j}
       +   \Quad^{2}_{ij}(X) \wedge \, X^{i}d X^{j}
\end{equation}
where $\Quad^1_{ij}$ and $\Quad^2_{ij}$ are smooth differential forms on
${\left.BM\right|}_{V} \subset TM$, which depend on the smooth form $\psi$
and on the coordinate chart $\phi_{\ell}$.
Because the maps $\Quad^a_{ij}$ are smooth differential forms
 for any smooth $q$-form
$\psi$, and because $M$ is compact,
we have the following  corollary to Lemma~\ref{lem:exp-q-form}, which we
prove in \cite{BD1}:

\begin{lemma}
\label{lem:remainder-estimate}
Let $\psi$ be a smooth $q$ form. Then the following estimates
are satisfied
for all $X \in \Wk{s}(TM)$, $s \geq 2n+6$, such that $|X|<c$:
\begin{subequations}
\begin{align}
\tag{a}
   \|\FX^*\psi \|_{s-2} &\prec
     \|\psi\|_{s-2} + \|\Lie_{X}\psi\|_{s-2} + \|X\|_{s-2} \, \|X\|_{s}
     \,,
\\
\tag{b}
   \|(\FX^*\psi) \wedge \n \|_{s-1} &\prec
     \|\psi \wedge \n\|_{s-1} + \|\Lie_{X}\psi \wedge \n\|_{s-1}
      + \|X\|_{s-1} \, \|X\|_{s} \,.
\end{align}
Moreover,  the estimate
\begin{equation}
\tag{c}
\begin{split}
   \|  \left(\Quad_{\psi} (X_1) - \Quad_{\psi} (X_2) \right)
\wg \n \|_{s-1} \prec&\quad
	 \|X_1 - X_2\|_{s-1} \, (\|X_1\|_{s} + \|X_2\|_{s})\\
&\quad +
	 \|X_1 - X_2\|_{s} \, (\|X_1\|_{s-1} + \|X_2\|_{s-1})\,,
\end{split}
\end{equation}
\end{subequations}
holds for any two  vector fields  $X_i$, $i=1,2$ with
$|X_{i}|<c$.
\end{lemma}

\begin{remark}
\label{rem:smooth-dependence}
As shown in \cite{BD1},  for $\psi$ a smooth $p$-form, the maps $X \mapsto
F^*_X\psi$  and $X \mapsto \n \wedge F^{*}_{X}\psi$ define smooth maps
\[
            \Wk{s}(TM) \to \Wk{s-2}(\Lambda^p M)\,,\text{ for  $s\geq 2n+6$},
\text{ and }
            \Wk{s}(TM) \to \Wk{s-1}(\Lambda^{p+1} M)\,,\text{ for $s \geq 2n+5$.}
\]
\end{remark}

Recall that  the condition for the diffeomorphism $\FX$ to be a
contact diffeomorphism is the vanishing of the one-form $\FX^* \n$ mod
$\n$. Hence by Equation~\eqref{eqn:calR-notation}, $\FX$ is a contact
diffeomorphism if and only if it satisfies the condition
\begin{equation}
\label{eq-contact-condition}
           \Lie_{X}\n + \Quad_{\n}(X) = 0  \text{ mod } \n \,.
\end{equation}
Furthermore, by Equation~\eqref{eq:Quad-local}, the linearization of this condition
at the zero vector field is the condition
\[
              \Lie_{X}\n = 0 \text{ mod } \n \,,
\]
i.e. $X$ is a \emph{contact vector field}.

\begin{remark}
\label{rem:contact-vf-param}
Using the characteristic vector field $T$ for the contact form $\n$, we may
express any vector field $X$ as $X = X^0 T + \XH $, where $\XH $ belongs to
the contact distribution.  Applying the Cartan identity
\[
	  \Lie_X (\n) = X \inter d \n + d (X \inter \n)
		      = \XH  \inter  d \n + d X^0 \,,
\]
yields the well known facts that (i) the vector field $X$ is a contact 
vector field if and only if
\begin{equation}
\label{contact.vf}
	d X^0 = - \XH  \inter d\n  \qquad \mod \n \, ;
\end{equation}
and (ii) that $X$ is completely determined by the real-valued function
$X^0 = X\inter \n$. For this reason $X\inter\n$ is called the \emph{generating function} for $X$
and is denoted by $g_X$. In \cite{BD1}, we proved that there is  an isomorphism
\[
         \Wkcont{s} \to \Wk{s+1}(M)\,:\, X \mapsto g_X := X\inter\n \,.
\]
\end{remark}
The main result of \cite{BD1} is the construction of a smooth
parameterization $\Psi$ of the space of $\Wk{s}$-contact
diffeomorphisms near the identity diffeomorphism by contact
vector fields near the zero vector field. The parameterization $\Psi$ in turn
induces a smooth structure on the space $\cdiff{s}$ of all 
$\Wk{s}$-contact diffeomorphisms.

\begin{theorem}[\cite{BD1}]
\label{thm:second-order}
For all $s \geq 2n+4$, and for $\cU \subset \Wk{2n+4}(TM)$ sufficiently small, there is a smooth map
\[
\Psi : \Wkcont{s}\cap \cU \rightarrow \cU^{s} \subset \Wk{s}(TM)
\]
such that the following holds: for all $Y \in \cU \cap \Wk{s}(TM)$, $F_{Y}$ is a contact
diffeomorphism if and only if $Y = \Psi(X)$ for some $X \in \Wkcont{s}
\cap \cU$.
Moreover, the map $\Psi$ is of the form
\[
	   \Psi(X) = X + B(X)(X,X) \, ,
\]
where $B: (\Wkcont{s} \cap \cU) \times \Wkcont{s} \times \Wkcont{s} 
	   \rightarrow \Wk{s}(TM)$ is smooth and bilinear in the last two factors.
\end{theorem}

This theorem implies  the following global result:
\begin{theorem}[\cite{BD1}]
\label{thm-4.20}
Let $(M,\n)$ be a compact contact manifold. For $s \geq (2n+4)$, the
space of $\Wk{s}$ contact diffeomorphisms is a smooth Hilbert manifold.
\end{theorem}
We close this section with the   \apriori\ estimates for the
nonlinear term $B(X)(X,X)$, which we proved in \cite{BD1} and which we require in Section~\ref{sec:slice}:

\begin{proposition}[\cite{BD1}]
\label{prop:Psi-estimates}
For $X \in \cV^{s}=\Wkcont{s} \cap \cU^s$,
\begin{equation}
\tag{a}
\label{sharper-estimate}
 \|\Psi(X) - X\|_{s} \prec \|X\|_{s}\|X\|_{s-1} \, .
\end{equation}
Moreover, for all  $ X_{1}, X_{2} \in   \cV^{s}$,
\begin{align}
\tag{b}
\label{cauchy-estimate}
\| (\Psi(X_{2}) - X_{2}) - (\Psi(X_{1})  - X_{1}) \|_{s}
\prec  &\; \|X_{2} - X_{1}\|_{s-1}(\|X_{2}\|_{s} + \|X_{1}\|_{s})\\
 &\; + \|X_{2} - X_{1}\|_{s}(\|X_{2}\|_{s-1} + \|X_{1}\|_{s-1}) \, .\nonumber
\end{align}
\end{proposition}

\section{The action of the contact diffeomorphism group}
\label{sec:action}

There is a natural action of contact diffeomorphisms on the space of
CR~deformations:
\[
   \begin{cases}
   \cdiff{\infty} \times  \FormsH{1} &\to  \FormsH{1}\\
      \qquad\quad     (F,\phi) &\mapsto    \mu = F^{\str}\phi
   \end{cases} \,.
\]
The main result of this section (Proposition~\ref{prop:pull-back-formula})
is a formula for $F^*\phi$  in the special case where $F = F_{\Psi(X)}$.

Let $F$ be a contact diffeomorphism, and let $\phi$ be a deformation
tensor.  Let $\Hhat_{(0,1)}\subset \HC = H_{(0,1)} \oplus H_{(1,0)}$ denote
the anti-holomorphic tangent bundle of the strongly pseudoconvex CR~structure associated to
$\phi$, and define the \emph{pull-back} CR~structure $F^*\Hhat_{(0,1)}
\subset \HC$ to be the CR~structure with anti-holomorphic subbundle
\[
      F^{*}(\Hhat_{(0,1)}) := \left\{ Z \in \HC \st   F_{*}Z \in \Hhat_{(0,1)}
      \right\} \,.
\]
It is straightforward to check that if $F_1$ and $F_2$ are two contact
diffeomorphisms then the identity
\[
            (F_2 \circ F_1)^{*} (\Hhat_{(0,1)}) = F_1^{*} (F_2^{*}(
            \Hhat_{(0,1)} ))
\]
holds. By Corollary~\ref{cor:deform-tensor}, if $F$ is isotopic to the identity,
then $F^{*}\Hhat_{(0,1)}$ is represented by a deformation tensor, which we
call the \emph{pull-back CR~deformation}, denoted by $F^{*}\phi$. 

\medskip

\subsection{Local formul\ae}
We need a local formula for $\FPX^{*}\phi$ that exhibits the non-linear
dependence on the contact vector field $X$.  It will also prove important to single out terms involving
composition of the components of the tensor $\phi$ with $\FPX$; we
accomplish this by introducing an auxiliary contact vector field $Y$ into
some \formulae.

Choose an adapted atlas and subordinate partition of unity as in
 Remark~\ref{rem:local-coords}.
By smoothness of the map $X \mapsto \FPX$ and compactness of $M$, for all
sufficiently small $X$, the condition
\[
      \FPX(V_{\ell}) \Subset U_{\ell}
\]
holds for  all $\ell$.
Next let $\n, \w_{\ell}^{\a}, \w_{\ell}^{\ab} = \overline
{\w_{\ell}^{\a}}$ be a  coframing for $T_{\C}M$ on $U_{\ell}$, with $\Hconj$ the
annihilator of $\n, \w^{\a}$.

For ease of notation, we temporarily suppress the index $\ell$ and set $F = \FPX$.
Then \begin{equation}
\label{eq:pullback-def}
      F^{*}(\Hhat_{(0,1)}) = \left\{ Z \in \HC \st   Z \inter  F^{*} (\w^{\a} +
        \phi^{\a}_{\bb} \w^{\bb}) = 0, \text{ for } \a = 1,2,\dots,n
      \right\}\,.
\end{equation}
One sees immediately that
\begin{equation}
\label{eqn:action-3}
      F^{*}(\w^{\a} + \phi^{\a}_{\bb}\w^{\bb}) 
	      =   A^{\a}_{\b} \w^{\b} + B^{\a}_{\bb}\w^{\bb}
\,\modeta \,,
\end{equation}
where
\begin{subequations}
\label{eqn:action}
\begin{eqnarray}
	\label{eqn:action-1}
	A^{\a}_{\b} &=& Z_{\b} \inter F^{*} (\w^{\a} +
    \phi^{\a}_{\gb} \w^{\gb} ) 
= \left( Z_{\b} \inter F^{*} \w^{\a} \right) +
    (\phi^{\a}_{\gb}\circ F )\, \left( Z_{\b}\inter F^{*} \w^{\gb} \right)
\\
	\label{eqn:action-2}  
    B^{\a}_{\bb} &=& Z_{\bb} \inter F^{*} (\w^{\a} +
    \phi^{\a}_{\gb} \w^{\gb} )
= \left( Z_{\bb} \inter F^{*} \w^{\a} \right) +
    (\phi^{\a}_{\gb}\circ F )\, \left( Z_{\bb}\inter F^{*} \w^{\gb} \right)
\,.
\end{eqnarray}
\end{subequations}
 By Lemma~\ref{lem:exp-q-form}, one has the \formulae\ 
\begin{subequations}
\begin{align}
\label{eqn:deform-0i}
    F^*\w^{\a}  &=  \w^{\a}  +  \Lie_{X}\w^{\a}  + \Lie_{(\Psi(X)-X)} w^{\a}
 + \Quad_{\w^{\a}} ( \Psi(X) ) \\
\intertext{and}
\label{eqn:deform-0ii}
   F^*\w^{\ab}  &= \w^{\ab} +   \Lie_{\Psi(X)}\w^{\ab} 
 +  \Quad_{\w^{\ab}} ( \Psi(X)) \,,
\end{align}
\end{subequations}
for one-forms $\Quad_{\w^{\a}}$ and $\Quad_{\w^{\ab}}$ as in Equation~\eqref{eq:Quad-local}
Consequently,
\begin{equation}
\label{eq:F*}
    F^*\left( \w^{\a} + \phi^{\a}_{\gb} \w^{\gb} \right)  = 
\w^{\a} + \Lie_{X}\w^{\a} + (\phi^{\a}_{\gb} \circ \FPX) \w^{\gb} + 
\Quad^{\a}(X,X,\phi) \,,
\end{equation}
where the expression $\Quad^{\a}(X,Y,\phi)$ is defined by the formula
\begin{equation}
\label{eqn:deform-A2}
\begin{split}
\Quad^{\a}(X,Y, \phi) := &  \Lie_{\Psi(X)-X}\w^{\a} +
    ( \phi^{\a}_{\gb}\circ \FPY)  \Lie_{\Psi(X)}\w^{\gb}\\
   &+  \left\{  \Quad_{\w^{\a}}(\Psi(X) ) + \phi^{\a}_{\gb}\circ \FPY)\, \Quad_{\w^{\gb}}  (\Psi(X))  
     \right\}
\end{split}
\end{equation}
for $Y$ a second, sufficiently small, contact vector field.

To single out the terms of the form $\phi^{\a}_{\gb}\circ \FPX$,
we replace  the term $\phi^{\a}_{\gb}\circ F$  in
Equations~\eqref{eqn:action-1} and \eqref{eqn:action-2} by
$\phi^{\a}_{\gb}\circ \FPY$ to get matrix-valued functions
\begin{equation}
\label{eqn:deformAB}
         A = A(X,Y,\phi) \text{ and } B = B(X,Y,\phi) \,.
\end{equation}
Using the identity
\[
Z_{\bb} \inter (\Lie_{X}\w^{\a}) = -
(\Lie_{X} Z_{\bb}) \inter \w^{\a} 
=  (\Lie_{Z_{\bb}} X) \inter \w^{\a} =
(\bDb X)^{\a}_{\bb} \,,
\]
and the expression for  $\bDb X$ in Remark~\ref{rem:lifted-operator},
 yields the following \formulae\  for the entries of $A(X,Y,\phi)$ and $B(X,Y,\phi)$
\begin{subequations}
\begin{align}
\label{eqn:deform-A1}
A^{\a}_{\b} &= \delta^{\a}_{\b} + Z_{\b}\inter\Lie_{X}\w^{\a}  + Z_{\b}\inter
         \Quad^{\a}(X,Y,\phi)
\\
\intertext{and}
\label{eqn:deform-B1}
B^{\a}_{\bb} &= (\bDb X)_{\bb}^{\a} + (\phi^{\a}_{\bb}\circ \FPY) + 
		  Z_{\bb} \inter \Quad^{\a}(X, Y,\phi) \,.
\end{align}
\end{subequations}
Finally, expressing $A^{-1}B$ 
in the form
$\displaystyle
	 A^{-1} B  =  B +  A^{-1} (I - A) B 
$
yields the identity
\[
    \FPX^*\phi =  \bDb X + \left( \phi_{\bb}^{\a}\circ \FPX  + \cE_{\bb}^{\a}(X,X,\phi)
    \right)\, \w^{\bb} \otimes Z_{\a}
\]
where
\begin{equation}
\label{action-6}
    \ \cE^{\a}_{\bb}(X, Y, \phi)  = Z_{\bb}\inter \Quad^{\a}(X,Y,\phi) +
     [A^{-1}(I-A)B]^{\a}_{\bb}\,.
\end{equation}
and where $A=A(X,Y,\phi)$, $B = B(X,Y,\phi)$.

Using the partition of unity, we globalize these local \formulae\ to obtain
the vector-valued one forms
\begin{equation}
\label{eqn:phi-circ-F}
	\phi\circ \FPX := \sum_{\ell}  
	\rho_{\ell} \cdot (\phi_{\ell,\bb}^{\a} \circ \FPX) \, 
	 \w^{\bb}_{\ell} \otimes  Z_{\ell,\a} 
\end{equation}
and
\begin{equation}
\label{eqn:EX-global}
	\cE(X,Y,\phi) := \sum_{\ell} \rho_{\ell} \cdot
        \cE^{\a}_{\ell,\bb}(X,Y,\phi)\,\w^{\bb}_{\ell}\otimes Z_{\ell,\a} \,.
\end{equation}
Noting that  $\mu = \sum_{\ell}\rho_{\ell} \cdot \mu$ and $\bdb X =
\sum_{\ell}\rho_{\ell} \bdb X$ then immediately gives the next proposition,
which we need to
prove the normal form theorem of Section~\ref{sec:slice}.

\begin{proposition}
\label{prop:pull-back-formula}
Let $(\FPX, \phi) \in \cdiff{s+1} \times \Wk{s}\left( (H^{(0,1)} \otimes
  \Hhol\right)$ be near $(id_{M},0)$.  Then $\FPX^{\str}\phi $
is given by
\begin{subequations}
\begin{equation}
\label{eqn:action-5}
\tag{a}
    \FPX^*\phi =  \bDb X + \phi\circ \FPX  + \cE(X,X,\phi) \,.
\end{equation}
The linearized action at the identity map and the zero deformation tensor
is
\begin{equation}
 (X, \dot \phi) \mapsto \bdb X + \dot \phi .
\label{eq:linearized-action}
\tag{b}
\end{equation}
\end{subequations}
\end{proposition}

\begin{remark}
\label{rem:non-invariance}
These equations require some care in interpretation. First, notice that the
terms $\FPX^{{\str}}\phi$ and $\bDb X$ are in fact globally defined tensors,
and make invariant sense. On the other hand, $\phi \circ \FPY$, and
$\cE(X,Y,\phi)$ have been defined using local coordinates and are
coordinate dependent.
\end{remark}

\begin{remark}
\label{rem:smoothing-of-cE}
Observe that for $s \geq
2n+4$ the map $(X,Y,\phi)\mapsto \cE(X,Y,\phi)$  extends  to the map
\[
     \cE:  \Wkcont{s+1} \times \Wkcont{s+1}\times  \Wk{s}\left( (H^{(0,1)} \otimes
            \Hhol\right) 
 \to   \Wk{s}\left( (H^{(0,1)} \otimes \Hhol\right) 
\]
between Folland-Stein spaces. In Section~\ref{sec:apriori} we
obtain estimates for $\cE$ that play a critical role in the proof or
our normal form theorem. 
\end{remark}

\subsection{Complex contact vector fields}

By equation~(\ref{eq:linearized-action}), the action of the group of
contact diffeomorphisms suggests normalizing deformation tensors by the
image of $\bdb X$ where $X$ is a real contact vector field.  On the other
hand, since $\bDb X = \bDb (\pi_{(1,0)} X)$, it is natural to normalize the
deformation tensor by the image of $\bdb X$ for $X \in \Thol$. We
accomplish this by
introducing the notion of complex
contact vector fields.\footnote{This corresponds to the notion of {\it
    Hamiltonian vector fields} as used in \cite{BD2}.}  

Begin by recalling that $\Thol$ is defined as the quotient bundle
\[
      0 \rightarrow \Hconj   \hookrightarrow \TMC
\stackrel{\pi_{(1,0)}}{\longrightarrow} \Thol \rightarrow 0\,.
\]
Noting that $\TMC = \Hhol \oplus \Hconj \oplus \C \cdot T$, where $T$ is the Reeb
vector field, we see that the restriction of $\pi_{(1,0)}$ to $\Hhol \oplus
\C \cdot T$ is an isomorphism of complex vector bundles. Thus, we shall
identify $\Thol$ with $\Hhol \oplus \C \cdot T$ when it is convenient.

Next observe that the composite map 
\[
      TM   \hookrightarrow \TMC
\stackrel{\pi_{(1,0)}}{\longrightarrow} \Thol
\]
is injective with image the subbundle $\{Z\in \Thol \st \n(Z) \in \R
\}$. Consequently, there are  natural identifications
\begin{equation}
\label{eqn:real-ident}
    TM = \Hhol \oplus \R \cdot T := \{ X \in \Thol \st   \n(X) \in \R \} \,;
\end{equation}
and it is easy to check that  the inclusion
$\Wk{s}(TM) \subset \Wk{s}(\Thol)$ is norm preserving.

Because $\Hconj$ is contained in the annihilator of $\n$, the quantity
$\n(Z)$ is well-defined for all $Z \in \Thol$. In addition, the quantity
$\pi^{(0,1)}(Z\inter d\n)$ is  well-defined, where $\pi^{(0,1)}:\TMC^*
\to H^{(0,1)}$ denotes the natural projection map. More precisely,
let\footnote{Independence of the choice of $\widetilde{Z}$ is an immediate consequence
of the identity $d\n(W_1,W_2) = 0$ for all $W_1,W_2 \in \Hconj$.}
\[
      \pi^{(0,1)} (Z \inter d\n) := \pi^{(0,1)} (\widetilde{Z} \inter
      d\n) \,,
\]
for any $\widetilde{Z} \in \TMC$ such that  $\pi_{(1,0)}\widetilde{Z} = Z$.

Finally, recall from Remark~\ref{contact.vf}
 that a real vector field $X$ is a contact vector field if and only
if it satisfies the identity
\[
       d X^0 + X \inter d\n = 0 \quad \mod \n \,,
\]
where $X^{0} = X\inter\n$. This is equivalent to the two conditions
\begin{quotation}
$W \inter \left(d X^0 + X \inter d\n \right) = 0 $  and 
$\bW \inter \left(d X^0 + X \inter d\n \right) = 0$ for all $W \in
\Hhol$.
\end{quotation}
Since $X$ is real, $W \inter \left(d X^0 + X \inter d\n \right)=\overline{W
  \inter \left(d X^0 + X \inter d\n \right)}$.  This leads us to the
following definition.

\begin{definition}
\label{def:complex-contact}
We say that a $(1,0)$-vector field $Z \in \Gamma(T_{(1,0)}M)$ is a
\emph{complex contact vector field} if it satisfies the condition
\[
    \bDb \left(Z \inter \n \right) + \pi^{(0,1)} (Z \inter d\n) = 0 \,.
\]
We denote by $\Wkcxcont{s}$ the Folland-Stein completion of the space of
complex contact vector fields.
\end{definition}

\noindent
The following lemma places this definition in context. 

\begin{lemma}
\label{lem:contact-char}
Let $X \in \G (T_{\C}M)$ be a real vector field. 
Then the following are equivalent:
\begin{enumerate}
\item[(a)] $X $ is a contact vector field.
\item[(b)] $X$ satisfies the identity
\[
    \bdb (X \inter \n) + \pi^{(0,1)}(X \inter d\n) = 0 \,
\]
where $\pi^{(0,1)}: T_{\C}^{*}M \rightarrow \Hconj^{*}M$ is the restriction
operator.
\item[(c)] The vector-valued one-form $\bDb X$ takes values in $\Hhol$.
\end{enumerate}
\end{lemma}

\begin{proof}
  By the observations above, a real vector field $X$ is contact if and only
  if
\[
\bW \inter \left(d (X \inter \n) + X \inter d\n \right) = 0 \qquad  \forall
\;
\bW \in \Hconj \; .
\]
This is simply condition \ref{lem:contact-char}(b), thus establishing the
equivalence of 
\ref{lem:contact-char}(a)
and
\ref{lem:contact-char}(b).
The equivalence of \ref{lem:contact-char}(b) and \ref{lem:contact-char}(c)
is a special case of Lemma~\ref{lem:complex-contact} below.
\end{proof}

The next lemma gives a  useful characterization of complex contact vector
fields.  Before stating the lemma, we remark that the quotient bundle $\Thol$ has
a naturally defined subbundle determined by the vanishing of $\n$, that
  is 
\[
            \HC / \Hconj = \{ Z \in \Thol \st \n(Z) = 0 \} \subset \Thol\,.
\]
A simple computation shows  that  the map $\pi_{(1,0)}$ defined above restricts to an
isomorphism  $\Hhol \simeq\HC/\Hconj$. Hence, we may identify 
$\HC/\Hconj$-valued forms with $\Hhol$-valued forms.

\begin{lemma}
\label{lem:complex-contact}
The vector field $Z \in \Gamma(\Thol)$ is a complex contact vector field if
and only if
 $\bDb Z$ is an \mbox{$\Hhol$-valued} $(0,1)$-form.
\end{lemma}
\begin{proof}
Suppose that $\bDb Z$ takes its values in $\Hhol$; that is, that 
\begin{equation}
	 \n (\bDb Z (\bW) ) = 0 \qquad \forall \; \bW \in \Hconj \; .
\label{eq:Hhol-valued}
\end{equation}
By equation~\eqref{eq:bDbX},
\[
  \n\left( \bDb Z (\overline{W}) \right) = 
-\n\left(\pi_{(1,0)}
[Z,\overline{W}] \right)
=
-\n\left( [Z,\overline{W}] \right) \; ;
\]
but
\begin{align*}
\n\left( [Z,\overline{W}] \right)
&= -d\n(Z,\overline{W}) + Z \n(\overline{W}) -
\overline{W}\n(Z)\\
&= -d\n(Z,\overline{W}) - \overline{W} (Z\inter\n)   \\
&= - \overline{W} \inter \left( Z\inter d\n  + \bDb
(Z\inter \n) \right)\,.
\end{align*}
Thus, $\bDb Z$ takes its values in $\Hhol$ if and only if $ \overline{W}
\inter \left( Z\inter d\n + \bDb (Z\inter \n) \right) = 0$ $\; \text{ for
  all }\bW \in \Hconj$, which is equivalent to $Z$ being complex contact.
\end{proof}


\section{Homotopy operators for CR manifolds}
\label{sec:homotopy}

In this section, we will collect various results concerning the existence
and regularity of homotopy operators on embedded strongly pseudoconvex CR
manifolds.  We restrict our statements to the special case of embedded,
three dimensional CR manifolds.   More details of these constructions and
their generalizations can be found in e.g.  \cite{BuMi}, \cite{Miy1}.

\subsection{Miyajima's homotopy operators}

First, we have the following result for the $\bdb$ complex. It follows
 immediately from the vector bundle valued version contained in
 \cite{Miy2}, where the vector bundle is the trivial line bundle, and $\Ph
 =
 N \bdb^{*}$, but the result is essentially contained in
 \cite{BeGr}. Roughly speaking, it says that there exists a partial
 inverse
 and a Szeg\"o projector with good estimates.

 \begin{theorem}
\label{thm:scalar-homotopy}
 There exist linear operators 
\[
\Hh : C^{\infty}(M) \rightarrow C^{\infty}(M)\,,\quad
\Ph: \Forms{1} \rightarrow  C^{\infty}(M)\,,
\text{ and } 
\Sh :  \Forms{1} \rightarrow \Forms{1} \,,
\]
such that the following identities and estimates are satisfied:
\begin{subequations}
\begin{align}
&\tag{a} 
   \bdb\circ \Hh = 0\,,\quad
   \Ph \circ \Sh = 0 \,,\quad 
   \Sh \circ \bdb = 0 \\
&\tag{b} 
  u = \Ph \bdb u + \Hh u \text{ and }
  \a =\bdb \Ph \a + \Sh \a \,,  \\
&\tag{c} 
 \|\Hh(u)\|_{s} \prec \|u\|_{s} \,\quad
  \|\Ph(\a)\|_{s} \prec \|\a\|_{s+1} \,,\quad
 \|\Sh(\a)\|_{s} \prec \|\a\|_{s} \,,
\intertext{for all $u \in C^{\infty}(M)$, $\a \in  \Forms{1}$, and $s\ge
  0$.}
&\tag{d} 
\text{ $\Hh$ extends to a self-adjoint, projection operator on $L^2(M,\C)$.}
\end{align}
\end{subequations}
\end{theorem}
\addtocounter{equation}{-1}  

Similarly, homotopy operators for $\Thol$-valued $(0,1)$ forms also exist,
with similar estimates \cite{Miy2}.  These estimates work in general for
vector valued forms, where the vector bundle is the restriction of a
complex vector bundle which extends to the complex manifold bounded by $M$
as a holomorphic bundle.  (If the complex space $X$ bounded by $M$ is
singular, we first resolve the singularities of $X$ and then apply the
above definition.)

\begin{theorem}[Miyajima \cite{Miy1}, \cite{Miy2}]
 \label{thm:Miyajima-3D}
There exist linear operators 
\begin{subequations}
\begin{align}
&\notag
  \rh : \Gamma^{\infty}(\Thol) \rightarrow
 \Gamma^{\infty} (\Thol)\,,\\
&\notag
  \Ph: \FormsT{1}  \rightarrow \Gamma^{\infty}(\Thol)\,,\\
&\notag
  \Qh : \FormsT{1} \rightarrow \FormsT{1}
\intertext{satisfying the following identities and estimates:}
&\tag{a}
\bDb \circ \rh = 0\,,\quad \Ph \circ \Qh = 0 \qquad \Qh \circ \bDb = 0\\
&\tag{b}
  Z = \Ph \bDb Z + \rh Z \text{ and } 
 \phi = \bDb P \phi + Q \phi\\
&\tag{c}
  \|\Ph \phi\|_{s+1} \prec \|\phi\|_{s} \,,
\qquad \|\Qh \phi\|_{s} \prec \|\phi\|_{s}\,,   
\qquad \|\rh(Z)\|_{s} \prec \|Z\|_{s} \,,\\
\intertext{for all $Z \in  \Gamma^{\infty}(\Thol)$, $\phi \in \FormsT{1}$,
and $s \ge 0$.}
\intertext{Finally, there exist linear operators}
&\notag
  \Lh: \FormsT{1}
 \rightarrow \Gamma^{\infty}(\Thol)
\text{ and }
\Nh : \Gamma^{\infty}(\Thol)   \rightarrow \Gamma^{\infty}(\Thol)\,,\\
\intertext{with $\Lh$  a smooth horizontal linear first order
differential operator such that}
\tag{d} 
   &\Ph = \Nh \circ \Lh \,,\\
\intertext{and $\Nh$ satisfies the estimate}
\tag{e}
   & \|\Nh(Z)\|_{s+2} \prec \|Z\|_{s}\,
\text{ for all $Z \in  \Gamma^{\infty}(\Thol)$, $s \ge 0$.}
\end{align}
\end{subequations}
\addtocounter{equation}{-1}  
\end{theorem}

\subsection{Homotopy operators for complex contact vector fields}
These homotopy \formulae\  do not single out contact vector fields in
any significant manner. We now show how to modify the homotopy operators in
order to do so. We begin by introducing the {\it
raising} and {\it lowering} operators induced by the nondegenerate two
form $d\n$:

\begin{definition}
\label{def:sharp}
The \emph{lowering operator} is the vector bundle
map
\[
     (\;)^{\flat}:  TM  \to H^{*} \st X \mapsto X^{\flat} = X \inter d\n
\]
whose restriction to $H\subset TM$ is an
isomorphism between the contact distribution and its dual space.
The \emph{raising operator} is the inverse 
\[
    (\;)^{\sharp}: H^{*} \to H \st \phi \mapsto \phi^{\sharp} \,.
\]
\end{definition}

\begin{remark}
\label{rem:complex-contact-vf2}
The maps   $(\;)^{\flat}$ and $(\;)^{\sharp}$ 
of Definition~\ref{def:sharp} induce (complex)  linear maps
\[
    (\;)^{\flat}: \Thol \to H^{(0,1)} \st Z_{(1,0)} \mapsto    Z_{(1,0)}
    \inter d\n
\]
and 
\[
      (\;)^{\sharp}: H^{(0,1)} \to  \Hhol \,,
\]
where the map $(\;)^{\sharp}$ is an isomorphism of complex vector
bundles. Observe that by construction, 
\begin{equation}
         Z = \n(Z)\, T +    (Z^{\flat})^{\sharp}
\end{equation}
for all $Z \in \Wk{\infty}(\Thol)$.  Notice also, that by
Definition~\ref{def:complex-contact}, $Z$ is an element of the space
$\Wkcxcont{\infty}$ of complex contact vector fields
if and only if it satisfies the identity
\begin{equation}
\label{def:complex-contact-two}
    \bDb \left(\n(Z) \right) + Z^{\flat} = 0 \,.
\end{equation}
Thus, every complex contact vector field is of the form
\begin{equation}
\label{eq:complex-decomposition}
        Z_f = f T - (\bDb f)^{\sharp}
\end{equation}
for $f$ a smooth complex valued function. Moreover, the
inclusion $TM \hookrightarrow \Thol$ induces  the inclusion
\[
           \Wkcont{\infty} \hookrightarrow \Wkcxcont{\infty}
\st X \mapsto Z_{g_X}\,,
\]
where $g_X = \n(X)$. (See Remark~\ref{rem:contact-vf-param}.)
\end{remark}

\begin{proposition}
\label{prop:proj-operators}
There exist smooth linear  operators
\[
     \hat{\cPh}, \hat{\cSh}\,:\, \Wk{\infty}(\Thol) \rightarrow
     \Wk{\infty}(\Thol)
\]
satisfying the following:
\begin{align}
&\tag{a}
\label{vf-decomp}
        Z = \hat{\cPh}Z + \hat{\cSh}Z\,
\text{ for all $Z \in \Wk{\infty}(\Thol)$.}\\
&\tag{b}
       \range(\hat{\cPh}) = \ker(\hat{\cSh}) = \Wkcxcont{\infty}\,.\\
&\tag{c}       
\label{eq:hat-identities}
    \hat{\cPh} \circ \hat{\cSh} = 0 \,,\quad
\hat{\cSh} \circ \hat{\cPh} = 0\,, \quad
	\hat{\cPh} \circ \hat{\cPh} = \hat{\cPh}\,, 
\text{ and } \hat{\cSh} \circ
	\hat{\cSh} = \hat{\cSh} \,.
\\
&\tag{d}   
 \|\hat{\cPh} Z\|_{s} \prec \|Z\|_{s} \text{ and }
 \|\hat{\cSh} Z\|_{s} \prec \|Z\|_{s} 
               \text{ for all $Z \in \Wk{\infty}(\Thol)$.}
\end{align}

\end{proposition}
\begin{proof}
Choose a vector field $Z \in \Gamma(\Thol)$, and compute as follows using
the homotopy operators from Theorems~\ref{thm:scalar-homotopy}
and~\ref{thm:Miyajima-3D}:
\begin{align}
\nonumber
	Z &= \n(Z)\,T + (Z^{\flat})^{\sharp}
           = (\Hh (\n(Z)) + \Ph(\bDb(\n(Z))))\,T 
                    + \left\{\bDb \Ph(Z^{\flat}) + \Sh(Z^{\flat})\right\}^{\sharp}\,.
                    \\
\intertext{Add and subtract the term $\Ph(Z^{\flat}) T$ and
  rearrange to get}
Z & =  \left\{( \Hh(\n(Z)) -\Ph(Z^{\flat}) )\,T + (\bDb
	\Ph(Z^{\flat}))^{\sharp} \right\}   
        + \left\{ \left( \Ph(\bDb( \n(Z))) + \Ph(Z^{\flat}) \right)T 
+ \left( \Sh(Z^{\flat}) \right)^{\sharp}\right\} \,.
\end{align}
Define  $\hat{\cPh}, \hat{\cSh}: \Wk{\infty}(\Thol)
\rightarrow \Wk{\infty}(\Thol)$ to be the linear operators given
by the \formulae
\begin{align*}
     \hat{\cPh}(Z) &= \left( \Hh(\n(Z)) -\Ph(Z^{\flat}) \right)\,T + \left(\bDb
	\Ph(Z^{\flat})\right)^{\sharp}\\
     \hat{\cSh}(Z) &= \left(\Ph(\bDb( \n(Z))) + \Ph(Z^{\flat}) \right)\,T 
+ \left( \Sh(Z^{\flat}) \right)^{\sharp} \,.
\end{align*}
By construction, $ Z =  \hat{\cPh} Z  +  \hat{\cSh} Z$ .

We claim that 
$\hat{\cP}Z$ is a smooth complex contact vector field. This follows from
Equation~\eqref{def:complex-contact-two} and the computation
\[
        \bDb\left( \n(\hat{\cPh}(Z))\right) +  \hat{\cPh}(Z)^{\flat}
        = \bDb\Hh(\n(Z)) - \bDb\Ph(Z^{\flat}) + \bDb\Ph(Z^{\flat}) =\bDb\Hh(\n(Z))=0\,.
\]
Observe also that by \eqref{def:complex-contact-two}
\[
\hat{\cSh}(Z) = \left( \Ph(\bDb( \n(Z))) + \Ph(Z^{\flat}) \right)\,T 
+ \left( \Sh(Z^{\flat}) \right)^{\sharp}
=  \left( \Ph( - Z^{\flat})  + \Ph(Z^{\flat})\right)\,T 
+ \left( \Sh( - \bDb (\n(Z))) \right)^{\sharp} = (-0)^{\sharp} = 0\,.
 \]
for all $Z \in \Wkcxcont{\infty}$.

We have shown that
  $\hat{\cPh}$ takes values in $\Wkcxcont{\infty}$ and that $\hat{\cSh}$
  vanishes on $\Wkcxcont{\infty}$. These facts, combined with
  Equation~\eqref{vf-decomp} imply that $\hat{\cPh}$ and $\hat{\cSh}$
  satisfy the identities :
\[  \hat{\cPh} \circ \hat{\cSh} = 0 \,,\quad
\hat{\cSh} \circ \hat{\cPh} = 0\,, \quad
	\hat{\cPh} \circ \hat{\cPh} = \hat{\cPh}\,, 
\text{ and } \hat{\cSh} \circ
	\hat{\cSh} = \hat{\cSh} \,,
\]
as well as the equalities
\[
       \Wkcxcont{\infty} = \range(\hat{\cPh}) =
\mathrm{ker}(\hat{\cSh}) \,.
\]
The estimates follow from the estimates in
Theorems~\ref{thm:scalar-homotopy} and 
\ref{thm:Miyajima-3D}.
\end{proof}

\begin{remark}
Because the projection operators $\hat{\cPh}, 
\hat{\cSh}$ in Proposition~\ref{prop:proj-operators}
 preserve the Folland-Stein regularity, they extend to projection operators on the
Folland-Stein space $\Wk{s}(\Thol)$ and they induce a direct sum decomposition
\[
\Wk{s}(\Thol) = \Wkcxcont{s} \oplus \ker(\hat{\cPh})
\]
with  $\Wkcxcont{s}= \ker(\hat{\cSh}) = \range(\hat{\cPh})  \subset
\Wk{s}(\Thol)$.
\end{remark}

The following variant of
Theorem~\ref{thm:Miyajima-3D}, highlights 
the role of contact vector fields.

\begin{theorem}
\label{thm:homotopy-3D}
 There exist linear operators 
\[
\cPh: \FormsT{1}  \rightarrow  \Wkcxcont{\infty}
\text{ and }
\cQh : \FormsT{1} \rightarrow \FormsT{1}
\]
such that:
\begin{subequations}
\begin{align}
&\tag{a} 
Z = \cPh \bDb Z + \rh Z \, \text{ for all $Z \in \Wkcxcont{\infty}$}\\
&\tag{b} 
     \phi = \bDb \cPh \phi + \cQh \phi\,
      \text{ for all  $\phi \in \FormsT{1}$}\\
&\tag{c}
     \bDb \cPh \circ \cQh = 0 \,, \qquad \cQh \circ \bDb \cPh = 0\\
&\tag{d}
     \cQh \circ \bDb {Z} = 0 \text{ for all $Z \in \Wkcxcont{\infty}$}\\
&\tag{e} 
 \|\cPh \phi\|_{s+1} \prec \|\phi\|_{s} \,,\quad \|\cQh
 \phi\|_{s} \prec \|\phi\|_{s}  \text{ for all  $\phi  \in \FormsT{1}$.}\\
\intertext{Moreover,}
&\tag{f} 
\|\rh(Z)\|_{s} \prec \|Z\|_{s} \, \text{ for all 
     $Z \in  \Wkcxcont{\infty}$,  $s \ge 0$.}
\intertext{Finally, there exists a smooth linear operators}
&\notag \cLh: \FormsT{1}
 \rightarrow \Gamma^{\infty}(\Thol)
\text{ and }
 \cNh : \Gamma^{\infty}(\Thol)   \rightarrow \Gamma^{\infty}(\Thol)
\\
\intertext{with $\cLh$ a horizontal, first order differential operator, such that}
&\tag{g}
   \cPh = \cNh  \circ \cLh\\
&\tag{h}
  \|\cNh(Z)\|_{s+2} \prec \|Z\|_{s} \,
\text{ for all  $Z \in  \Gamma^{\infty}(\Thol)$,  $s \ge 0$.}
\end{align}
\end{subequations}
\end{theorem}
\addtocounter{equation}{-1}  

\begin{proof}
The key step in the proof is to express the homotopy operator $P$ of
Theorem~\ref{thm:Miyajima-3D} as the sum of two operators $\cPh$ and
$\cSh$,  defined by the formulas
\[
            \cPh = \hat{\cPh} \circ \Ph \text{ and }
            \cSh = \hat{\cSh} \circ \Ph \,.
\]
By Proposition~\ref{prop:proj-operators}, $P = \cPh + \cSh$ and 
 the image of $\cP$ is contained in the space
$\Wkcxcont{\infty}$ of smooth complex contact vector fields. 
Next let $\cH = \bDb \circ \cS + \Qh$, where $\Qh$ is as in
Theorem~\ref{thm:Miyajima-3D}. 

To prove (a), let $Z$ be a complex contact vector
field and note that by \ref{thm:Miyajima-3D}(b)
\[
       Z = \Ph \bDb Z + \rh Z = \cPh \bDb Z + \cSh \bDb Z + \rh Z \,.
\]
We need only show that $\cSh \bDb Z = 0$, for $Z$  complex contact.
First observe that whenever $Z$ is a
 complex contact vector field, then $\Ph \bDb Z$ is also complex contact.
 This follows easily from Lemma~\ref{lem:complex-contact}, the formula $\Ph
 \bDb Z = Z - \rh(Z) $, and $\bDb \rh(Z) = 0$.  Consequently, $\cSh(\bDb Z)
 = \hat{\cSh}(\Ph \bdb Z) =0$, for all $Z \in \Wkcxcont{\infty}$.

To prove the homotopy formula (b), notice that
Proposition~\ref{prop:proj-operators}(a) implies the decomposition
\[
           \Ph = \cPh + \cSh \,;
\]
then use the homotopy formula \ref{thm:Miyajima-3D}(b) to compute as follows:
\[
 \phi = \bDb P \phi + \Qh \phi = \bDb \cPh + \bDb \cSh \phi + \Qh\phi
               = \bDb \cPh \phi + \cH\phi\,.
\]

We now prove parts (c) and (d).  
First observe that $\cSh \circ \bDb \cPh = 0$.
Since $\cPh\phi$ is   complex contact,  $\Ph\bDb\cPh\phi$ is   
complex contact. Therefore, $\cSh\bDb\cPh\phi =  \hat{\cSh} (\Ph\bDb\cPh\phi)= 0$.
Next observe  that  $\bDb (\cPh \circ \bDb \cSh) = 0$ as follows:
For $\phi \in \FormsT{1}$, compute as follows:
\[
	\bDb\cPh \phi  =
\bDb\cPh \left( \bDb \cPh \phi  + \bDb \cSh \phi +  \Qh \phi      \right)
   = \bDb\cPh  \bDb \cPh \phi  + \bDb\cPh \bDb \cSh \phi \, ;
\]
on the other hand
\[
	\bDb\cPh \phi  =  \bDb\cPh(\bDb\cPh\phi) + \bDb\cS(\bDb\cPh\phi)
	+  \Qh  (\bDb\cPh\phi) =  \bDb\cPh \bDb\cPh\phi   \,.
\]
Thus, $\bDb (\cPh \circ \bDb \cSh) = 0$.
Finally,  the identities $\Qh
\circ \bDb = \cPh \circ \Qh = \cSh \circ \Qh = 0$ follow immediately from
Proposition~\ref{prop:proj-operators} and Theorem~\ref{thm:Miyajima-3D}.
Then the identities $\bDb\cPh
\circ \cQh =0$ and $ \cQh \circ \bDb\cPh = 0$ follow
from the identities $\bDb (\cPh \circ \bDb \cSh) =0$ and $\cSh
\circ \bDb \cPh = 0$.

To prove part (g), set $\cLh = \Lh$ and $\cNh =
\hat{\cPh} \circ \Nh$.  Since $\cPh = \hat{\cPh} \circ \Ph$ and by
(\ref{thm:Miyajima-3D}c) $\Ph= \Nh \circ \Lh$, it follows that $\cPh = \cNh
\circ \cLh$. 

The estimates (e), (f), and (h) follow immediately from the estimates in
Theorems~\ref{thm:scalar-homotopy} and~\ref{thm:Miyajima-3D}.
\end{proof}

Notice that in the last theorem, since $\FormsH{1} \subset \FormsT{1}$, it
follows that for $\phi \in \FormsH{1}$, we have $\phi = \bDb \cPh \phi +
\cQh \phi$.  Moreover, since the range of $\cPh$ is the space of
complex contact vector fields, then $\bDb \cPh \phi \in \FormsH{1}$
(see Lemma~\ref{lem:complex-contact}).  It follows that $\cQh$ restricts to
an operator $\cQh : \FormsH{1} \rightarrow \FormsH{1}$. Therefore, we can
restrict the homotopy formula to the horizontal vector valued forms. We
state this next, using the same symbols to denote the restricted operators
without risk of confusion.

\begin{corollary}
\label{cor:homotopy-contact}
 There exist homotopy operators
 $\cPh: \FormsH{1}  \rightarrow \Wkcxcont{\infty}$,
\newline 
    $\cQh : \FormsH{1} \rightarrow \FormsH{1}$ 
such that:
\begin{subequations}
\begin{align}
  \tag{a} 
       &\phi = \bDb \cPh \phi + \cQh \phi 
           \text { for all $\phi \in \FormsH{1}$}\\
 \tag{b} 
     &\bDb \cPh \circ \cQh = 0 \qquad \cQh \circ \bDb \cPh = 0\\
 \tag{c} 
     &\cQh \circ \bDb {Z} = 0\,
             \text{ for all  $Z\in \Wkcxcont{\infty}$}\\
 \tag{d} 
     &\|\cPh \phi\|_{s+1} \prec \|\phi\|_{s} \qquad \|\cQh
           \phi\|_{s} \prec \|\phi\|_{s}  
        \text{ for all $\phi  \in \FormsH{1}$.}\\
\intertext{Moreover, noting that the harmonic projection $\rh$ restricts to a
 map $\rh : \Wkcxcont{\infty} \rightarrow \Wkcxcont{\infty}$:}
 \tag {e}  
    &Z = \cPh \bDb Z + \rh Z 
              \text{ for all $Z \in  \Wkcxcont{\infty}$}\\
 \tag {f}  
    &\|\rh(Z)\|_{s} \prec \|Z\|_{s}
 \text{ for all  $Z \in  \Wkcxcont{\infty},  s \ge 0$.}
\end{align}
\end{subequations}
\addtocounter{equation}{-1}  
\end{corollary}

\subsection{Harmonic decomposition of complex contact vector fields}
\label{subsect:splitting-cont-vf}

In this section, we obtain a decomposition of complex contact vector fields
into real contact vector fields and a complementary subspace.
Recall from Equation~\eqref{eq:complex-decomposition} that the space of
complex contact vector fields is parameterized by complex valued functions
as follows:
\[
     f \mapsto  Z_f = f T - (\bDb f)^{\sharp} \,.
\]

The observation that this  parametrization agrees with the parametrization of
real contact vector fields as introduced in
Remark~\ref{rem:contact-vf-param}  suggests constucting the
decomposition using the na\"{i}ve projection operator $\pi_{\Re}: Z_{f}
\mapsto Z_{\Re(f)}$.
Unfortunately, this projection map is not continuous in the Folland-Stein
norm.  We see this as follows. By virtue of the identification $\Thol =
\Hhol \oplus \C \cdot T$, the Folland-Stein structure on the space of
complex contact vector fields is
\begin{eqnarray*}
\label{}
(Z_{f}, Z_{g})_{s}  \; &= & \left( (f T - (\bDb f)^{\#}), ( g T - (\bDb
g)^{\#}) \right)_{s} \\
&= & \left( f T ,  g T \right)_{s} + \left(  (\bDb f)^{\#},  (\bDb
g)^{\#} \right)_{s} \\
&= & \left( f ,  g \right)_{s} + \left(  \bdb f ,  \bdb
g \right)_{s} .\\
\| Z_{f}\|^2_{s}  \; &= &  \| f \|^2_{s} + \| \bdb f \|^2_{s} \, .
\end{eqnarray*}
On the other hand, since  $Z_{\Re(f)} = 1/2 (Z_{f} + Z_{\overline{f}})$,
\[
       \|Z_{\Re(f)}  \|^2_{s} = 1/2 \|Z_{f} + Z_{\overline{f}} \|^2_{s}
                   = \|\Re(f) \|^2_{s} + 1/2 \| \bDb f +
                   \overline{ \partial_b f} \|_{s}^2 \,.
\] 
Let $f_k$ be a sequence of CR~functions with $\|\partial
f_k\|_{s}\to\infty$ and $\norm{f_k}_s$ bounded.
Then $\|Z_{f_k}\|_s$ is bounded, but 
\[
\|Z_{\Re(f_k)}\|^2_s = \|\Re(f_k)\|^2_s + \frac{1}{2} \|\bDb \overline{f_k}\|^2_{s} 
 = \|\Re(f_k)\|^2_s + \frac{1}{2}\|\partial_b f_k\|^2_{s} \to \infty\,.
\]
Therefore, to obtain a bounded projection, we have to proceed
differently. We need the following regularly lemma.
\begin{lemma}
\label{lem:reg-real}
 The estimate $\norm{u}_{s+2} \prec \norm{\Re(I + \Box_b) u}_{s}$ holds for
 any smooth, real-valued function $u$. In particular, if $\Re(I + \Box_b) u$
 is smooth, then so is $u$.
\end{lemma}
\begin{proof}
One easily verifies that for $u$ real, the identity
$\Re ( u + \Box_b u) =  u + \frac{1}{2n+2} \DeltaQ u$ holds, where $\DeltaQ$ is the
Laplace operator in the Rumin complex. The estimate follows from the
corresponding estimate for $\DeltaQ$, proved in \cite{Rumin:1994, BD3}.
\end{proof}

Next let $f$ be  a smooth,
complex valued function $f$. Then $\Re( f + \Box_b f)$ is smooth, and
Lemma~\ref{lem:reg-real} implies that there is a unique, smooth, real-valued function
$u$, satisfying the equation
\[
           (I + \frac{1}{2n+2}\DeltaQ) u =   \Re(f + \Box_b f) \,.
\]
\begin{proposition}
\label{prop:piRe}
For all $s \geq 2n+4$,  the map 
\[
         f \mapsto u := (I + \frac{1}{2n+2}\DeltaQ)^{-1}\Re(f + \Box_b f)
\]
induces a bounded projection operator
\[
        \pi_{\Re}: \Wkcxcont{s} \to \Wkcxcont{s} \st Z_{f} \mapsto Z_{u}
\]
with image $\Wkcont{s}$.
\end{proposition}
\begin{proof}
By construction $\pi_{\Re}(Z_u) = Z_u$ for $u$ real. Consequently,
$\pi_{\Re}$ is a projection operator, as claimed. 
To prove that $\pi_{\Re}$ is bounded, note that regularity
  for $\DeltaQ$ justifies estimating as follows:
\[
     \norm{Z_u}_{s} \prec \norm{u}_s + \norm{\bdb u}_s \prec
\norm{u}_{s+1} \prec \norm{ (I + \frac{1}{2n+2}\DeltaQ)u}_{s-1}
\prec \norm{ \Re( f + \Box_b f)}_{s-1}
\prec \norm{f}_{s-1} + \norm{\Box_b f}_{s-1} \,.
\]
But
$\norm{\Box_b f}_{s-1} = \norm{\bdb^* \bdb f}_{s-1} 
\prec \norm{\bdb f}_{s}$ implies the estimate
$ \norm{Z_u}_{s} \prec \norm{f}_{s-1} + \norm{\bdb f}_{s} \prec
          \norm{Z_{f}}_{s}$.
\end{proof}

The projection map $\pi_{\Re}$ induces the decomposition
\[
      \Wkcxcont{\infty} = \Wkcont{\infty} \oplus i V \,,
\]
where
\begin{equation}
\label{eqn:def-V}
    V : = \{ Y\in \Wkcxcont{\infty} \st \ \pi_{\Re} (i Y) = 0\} \,. 
\end{equation}
Let $V^s$ denote the closure of $V$ in the $\Wk{s}$ norm.
It will prove convenient to adopt the notational convention
\begin{equation}
\label{eq:complex-contact-decompose}
	   Z_{f} = X_{f} - i \, Y_{f} \,,
\end{equation}
where $X_{f} := \pi_{\Re}(Z_f) \in  \Wkcont{\infty}$ and  $Y_{f} :=
\pi_{\Im}(Z_f)$ is the projection
\[
   \pi_{\Im} := i(\Id   - \pi_{\Re}) : Z_{f} \mapsto Y_{f} \; .
\] 
Moreover, the estimate
\[
       \|X_f \|_s + \|Y_f \|_s \prec \|Z\|_{s}
\]
holds for all $f \in \Wk{s}(M,\C)$, with $s \geq 2n + 4$.

\begin{remark}
  We caution the reader that although $X_{f}$ is real, it is \underline{not} the
  real part of $Z_{f}$.
\end{remark}

\begin{remark}
\label{prop_of_V}
We could at this point let $iV$ be a rather arbitrary complement to 
 $ \Wkcont{\infty}$. The only properties for $V$ that are important in what
 follows  are:
\begin{enumerate}
  \item[(a)] $ \Wkcxcont{\infty} \cong \Wkcont{\infty}  \oplus iV$,
  \item[(b)] $\Wkcxcont{s} \cong \Wkcont{s}  \oplus iV^s$,
  \item[(c)] $\|X\|_s + \|Y\|_s \prec \|X - i Y\|_{s}$,
\end{enumerate}
 for  all $s \ge  2n+4$.
\end{remark}


\section{Normal form for CR~deformations}
\label{sec:slice}

In this section, we study the action of the contact diffeomorphism group on
the space of deformations of a fixed embeddable CR~structure $(M, \Hhol)$
on a compact three dimensional manifold $M$.

There are significant differences in the analysis between the three
dimensional case and higher dimensions. These arise since first, there are
no integrability conditions in dimension three, and second, the relevant
operators are not subelliptic in three dimensions.  While the analysis
generalizes to higher dimensions, the details are numerous and everything
requires a separate statement, including the introduction of new operators
to take into account the integrability conditions. Since in dimensions at
least five, it is well known that all compact, strongly pseudoconvex CR
manifolds are embeddable, our main interest is in the three dimensional
case where the situation is more subtle and less well
understood. Henceforth, we will restrict our attention to this case.

Before beginning the statement and proof of the main results, we make some
comments to motivate the definitions and statements. The contact
diffeomorphism group acts on the space of deformation tensors, and the
linearization of the action at the identity map and the zero deformation
tensor is $(X, \dot \phi) \mapsto (\bdb X + \dot \phi)$, where $X$ is a
contact vector field. On the other hand, the Hodge decomposition of
Corollary~\ref{cor:homotopy-contact} shows that a deformation tensor can be
split as $\phi = \bDb \cPh \phi + \cQh \phi$, where $\cPh \phi $ is a
complex contact vector field, and $\cQh \phi$ serves as the ``harmonic
part'' of the deformation. If we split the complex contact vector fields as $\cPh \phi = X - i
Y$, where $X$ is a real contact vector field, and $Y$ lies in a transverse
subspace (see Section~\ref{subsect:splitting-cont-vf}), then $Y$ can be
heuristically thought of as infinitesimally arising from one of Kuranishi's
``wiggles'' of the embedded CR manifold within its ambient surface. 
The normal form should then be $i \bdb Y + \phi_{\cQh}$, that is, a harmonic
form plus a wiggle.

This overview suggests that we should consider a map $\Wkcont{\infty}
\oplus iV \oplus \ker\cPh \rightarrow \Deform$ and show that for all $\phi \in
\Deform$, there exist $(X, Y, \psi)$ such that $\FPX^{*} \phi = i \bdb Y +
\psi$ or $(\FPX^{-1})^{*}(i \bdb Y + \psi) = \phi$; here, $\FPX$ is the
contact diffeomorphism defined by $\Psi(X)$ as in
Theorem~\ref{thm:second-order}.  Unfortunately, the linearization of this
map loses regularity, since it involves differentiation with respect to
$X$, which has a component in the direction of $X$.  To circumvent this
difficulty, we carry along a copy of $\phi$ and consider the modified map
$(\phi, X, Y, \psi) \mapsto (\phi, \FPX^{*} \phi - (i \bdb Y +
\psi))$. This map is now invertible (modulo a kernel -- the CR vector
fields -- which is easily incorporated) giving a weak normal form: 
\begin{quotation}
\em for
every $\phi$, there is a triple $(X, Y, \psi)$ such that $\FPX^{*} \phi = i
\bdb Y + \psi$. 
\end{quotation}
However, in the proof, the normal form $i \bdb Y + \psi$
has less regularity than $\phi$. This can be viewed as a weak Hodge
decomposition for the nonlinear theory.  However, taking our lead from the
proof of regularity for the standard linear Hodge theory, we then obtain
\apriori\ estimates in Section~\ref{sec:apriori} to improve the regularity
and establish a strong normal form:
\begin{quotation}
\em  if $\FPX^{*} \phi = i \bdb Y + \psi$
with $\phi \in \Wk{s}$, then $X, Y \in \Wk{s+1}, \psi \in \Wk{s}$.
\end{quotation}

\noindent
\textbf{Remark.} We expect that this approach of first using 
linear analysis to obtain a weak normal form and then \apriori\ 
estimates to obtain the strong normal form will find a wide range of 
use in other applications.

\subsection{Statement of the Normal Form Theorem}
\label{subsect:normal-forms}

Throughout the remainder of the paper, $(M, \Hhol)$ is a fixed embeddable
compact three dimensional CR manifold.

We first establish notation.  Let $\Harm = \ker \cPh \subset \FormsH{1} $
represent the ``harmonic deformation tensors'', where $\cPh$ is as in
Corollary~\ref{cor:homotopy-contact} and denote the the CR vector fields by $\WkCR{s+1} = \ker\bDb \cap
\Wk{s+1}(\Thol)$. Let $\Wk{s}(\Harm)$ denote the
Folland-Stein completion of $\Harm$ in $\Wk{s}( \FormsH{1} )$. Notice that
$\Wk{s}(\Harm) $ is closed in the space of deformation tensors
$\Wk{s}(\Deform) =\Wk{s}( \FormsH{1} )$ and that by
Corollary~\ref{cor:homotopy-contact}
\[
  \Wk{s}( \FormsH{1} ) = \range(\bDb) \oplus \Wk{s}(\Harm) \,.
\]
We define the map\footnote{\noindent This corrects a mistake in \cite{B}
  when we mistakenly asserted the map $\Phi$ to be $C^1$ if we take the
  first factor on each side to be in $\Wk{s}(\Deform)$. In the
  Section~\ref{sec:apriori}, we obtain \apriori\ estimates to establish a
  local nonlinear Hodge theory and recover the lost regularity.}:
\begin{subequations}
\begin{equation}
\label{eqn:Phi-defn-a}
   \Phi : \Wk{s+2}(\Deform) \oplus  \Wkcont{s+1} \oplus V^{s+1} \oplus 
   \Wk{s}(\Harm) 
 \longrightarrow 
\Wk{s+2}(\Deform) \oplus  \Wk{s}(\Deform) \oplus  \WkCR{s+1} 
\end{equation}
by the formula
\begin{equation}
\label{eqn:Phi-defn-b}
      (\phi, X, Y, \psi) \mapsto 
\left( \phi, \FPX^{{\str}}\phi - i\bDb Y - \psi, \rh (X-iY)  \right).
\end{equation}
\end{subequations}

\begin{proposition}
\label{prop:Phi-local-diffeo}
The map $\Phi$ is a local diffeomorphism in a neighbourhood of the
origin.
\end{proposition}

\begin{proof}
By the inverse function theorem for Banach spaces, it is sufficient to
establish that:
\begin{quotation}
\item{} (1) $\Phi$ is locally $C^{1}$;
\item{} (2) $d\Phi|_{(0,0,0,0)}$ is invertible.
\end{quotation}
  To establish (1), notice that all terms in the map $\Phi$ are linear,
  and  smooth (see Theorem~\ref{thm:Miyajima-3D}), except $\FPX^{{\str}}\phi$, so it suffices to check the
  regularity of this term.  By Remark~\ref{rem:smooth-dependence},  $X
  \mapsto (Z_{\b} \inter \FPX^{*}\w)$ and $X \mapsto (Z_{\b} \inter
  \FPX^{*}\bw)$ are smooth maps from $\Gamma^{s+1}$ contact vector fields
  to $\Gamma^{s}$ functions.  We proved in \cite{BD1} that the
  map\footnote{\noindent Notice that the differential of this map is the
    Lie derivative of $u$, which explains the loss in regularity on $u$.}
\[
\Gamma^{s+2}(M) \oplus \cdiff{s+1} \rightarrow \Gamma^{s}(M)
\st (u, F) \mapsto u \circ F
\]
is $C^{1}$.  From the local expressions in \formulae
~\eqref{eqn:action-1}
and \eqref{eqn:action-2} and  the fact that the matrix $A$ in these \formulae\ is invertible, it
follows that the term $(\phi, X) \mapsto \FPX^{{\str}}\phi$ is $C^{1}$,
completing the proof that the map $\Phi$ is $C^{1}$.

We next check that $d\Phi$ is invertible at the origin. Let
$(\dot\phi, \dot X, \dot Y, \dot\psi )$ be a tangent vector
at the origin. Then
\[
    d\Phi ( \dot\phi, \dot X, \dot Y, \dot\psi) = \left(\dot\phi ,
    \bDb \dot X - i \bDb \dot Y + \dot\phi - \dot\psi,
    \rh (\dot X - i \dot Y) \right)   \, .
\]
It is clear that this map has trivial kernel and that it is surjective.
In fact, using the homotopy operators $\cPh, \cQh$, we can verify that
the
 inverse map $(d\Phi)^{-1}$ is given by
\[
(d\Phi)^{-1}: \begin{cases}
&\Wk{s+2}(\Deform) \oplus \Wk{s}(\Deform) \oplus \WkCR{s+1}
 \longrightarrow
\Wk{s+2}(\Deform) \oplus \Wkcont{s+1} \oplus
V^{s+1} \oplus  \Wk{s}(\Harm) \\
&(\dot\phi, \chi, \xi)  
 \mapsto  (\dot\phi,
        {\pi_{\Re}}(\cPh(\chi -
	\dot\phi) + \xi),
	{\pi_{\Im}}(\cPh(\chi - \dot\phi)+ \xi),
	- \cQh(\chi -     \dot\phi))  
\end{cases} \,.
\]
To verify that this  is the inverse of $d\Phi_{(0,0,0)}$, compute as
follows:
\begin{align*}
(d\Phi)^{-1}(\dot\phi,\bDb \dot{X} - i \bDb \dot Y + \dot\phi - \dot\psi,
    \rh( \dot X - i \dot Y ) )   
 = &\;
	(\dot\phi, {\pi_{\Re}}(\cPh( \bDb \dot{X} - i \bDb \dot Y
	 - \dot\psi) + \rh(\dot X - i \dot Y ) ),\\
 & \;\; {\pi_{\Im}}(\cPh( \bDb \dot{X} - i \bDb \dot Y - \dot\psi)
	+ \rh( \dot X - i \dot Y )),\\ 
  & \;\; - \cQh( \bDb \dot{X}- i \bDb \dot{Y} - \dot\psi)) \\
 = &\; (\dot\phi, \dot X, \dot Y, \dot\psi) \,.
\end{align*}
\end{proof}

By the implicit function theorem, inverting $\Phi$ gives rise to the $C^{1}$ map
\begin{subequations}
\label{eq:Phi-inverse}
\begin{equation}
\label{eq:Phi-inverse-a}
  \Wk{s+2}(\Deform) 
   \rightarrow   \cdiff{s+1}  \oplus V^{s+1} \oplus
   \Wk{s}(\Harm)
\,:\,   \phi \mapsto   \left( F_{\phi}, Y_{\phi}, \psi_{\phi} \right)
\end{equation}
defined by the constraint
\begin{equation}
\label{eq:Phi-inverse-b}
(\phi, X_{\phi}, Y_{\phi}, \psi_{\phi})  =\Phi^{-1}(\phi, 0, 0)\; ,
\end{equation}
with $F_{\phi} = F_{\Psi X_{\phi}}$ and $\phi$ in a sufficiently small neighbourhood of the origin.
\end{subequations}

\begin{corollary}
\label{cor:weak-normal-form}
There exist neighbourhoods $0\in U \subset \Wk{s+2}(\Deform)$ and
$id_{M}\in \tilde U \subset \cdiff{s+1}$ such that for any $\phi \in U$,
there is a contact diffeomorphism $F_{\phi} \in \widetilde{U}$ such that $F_{\phi}^{\str}\phi$ is
contained in the subspace $\bDb\left(i V^{s+1}\right) \oplus \Wk{s}(\Harm)
\subset \Wk{s}(\Deform)$. The equation
\[
F_{\phi}^{{\str}} \phi =  i\bDb Y_{\phi} + \psi_{\phi} \in
\Wk{s}(\Deform)
\]
determines $F_{\phi}, Y_{\phi}$ and $\psi_{\phi}$ up to the $CR$-vector
field $\rh(X_{\phi} - i Y_{\phi})$, which is in turn determined by the additional constraint
$\rh(X_{\phi} - i Y_{\phi}) = 0$.
\end{corollary}

We call the deformation tensor 
\[
	    F_{\phi}{\str}\phi = i \bDb Y_{\phi} + \psi_{\phi} \in
            \Wk{s}(\Deform)
\]
the \emph{normal form} of $\phi$.  The following theorem, which is proved
using \apriori\ estimates, gives increased regularity for the normal form.
It is an immediate corollary to Theorem~\ref{apriori-estimates-1} below.

\begin{theorem}
\label{strong-normal-form}
The map $\phi \mapsto  (F_{\phi}, Y_{\phi}, \psi_{\phi})$  defines a
$C^{0}$ map of the form
\[
   \Wk{s+2}(\Deform)
   \rightarrow  \cdiff{s+3}
    \oplus V^{s+3} \oplus   \Wk{s+2}(\Harm)\,,
\]
for sufficiently small $\phi \in \Wk{s+2}(\Deform)$.  In particular, the
normal form
\[
F_{\phi}^{{\str}} \phi =  \left( i\bDb Y_{\phi} + \psi_{\phi} \right)
\]
is contained in $\bDb\left(i V^{s+3}\right) \oplus \Wk{s+2}(\Harm)
\subset \Wk{s+2}(\Deform)$.
\end{theorem}

\begin{remark}
  As noted in Remark~\ref{prop_of_V}, we have some freedom in the choice of
  $V^{s+3}$, the complementary subpace to $\Wkcont{s+3}$ in
  $\Wkcxcont{s+3}$. If the original CR manifold admits a free $S^1$ action
as a symmetry, we can choose all homotopy operators to be $S^1$
equivariant. Complex contact vector fields then have Fourier expansions,
and we can choose our complement $V$ to consist of complex vector fields of
the form $Z_f$, where $f$ has only positive (respectively negative) Fourier coefficients. In
\cite{B}, we made these choices to obtain the interior (respectively
exterior) normal form.

  In general, since $M$ is embeddable it follows that $M \hookrightarrow
  \Sigma$ for some compact complex surface $\Sigma$ as a separating
  hypersurface (see \cite{Lempert:alg-approx}.)  The elements of $V$
  correspond on the infinitesimal level to Kuranishi's ``wiggles'', that
  is, CR structures which are induced on $M$ through infinitesimal
  isotopies of $M$ within $\Sigma$. In this regard, one expects the factor
  $\psi$ to correspond to deformations of the singularities of the
  ``fill-in'' of $M$ (that is, the pseudoconvex side of $\Sigma$ bounded by
  $M$) or to non-embeddable structures on $M$.
\end{remark}


\subsection{A priori estimates for the action on CR~structures}
\label{sec:apriori}

We now proceed to establish the \apriori\ regularity estimates for
the action of the contact diffeomorphism group on the space of
deformation tensors that we need to
establish Theorem~\ref{strong-normal-form}.

Let $X$ be a contact vector field and let $\phi$ be a CR~deformation,
expressed relative to a local frame $Z_{\a}$ and
dual coframe $\w^{\bb}$ as $\phi =
\phi^{\a}_{\bb} \w^{\bb}\otimes Z_{\a}$.  For $X$ and
$\phi$ sufficiently small, we will obtain estimates for the 
deformation tensor  for the pull-back CR-structure $\mu = F^{\str}\phi$.
\footnote{Although we have restricted to the three dimensional case $n=1$, we
  continue to use index notation to help distinguish between functions and
  coefficients of tensors.}

\begin{remark}
\label{rem:neighbourhood}

Since we are restricting ourselves to a small neighbourhood of the
embeddable structure, we may choose the neighbourhood small enough to
have the following uniform estimates:
\begin{equation*}
  \|\phi\|_{s+2}   <  \const{\phi}{s+2} \,,
\quad
  \|\mu\|_{s+2} <   \const{\mu}{s+2}   \,,
\text{ and }
  \|X\|_{s+1} <  \const{X}{s+1}
\end{equation*}
where $\const{\phi}{k}$ is a fixed (sufficiently small) constant.  Since
when $ \|X\|_{s+1} < \const{X}{s+1}$ one has $ \|X\|_{s+1} \sim
\|\FPX\|_{s+1} $, and one can choose $C$ such that in addition
\[
    \|\FPX\|_{s+1}  <  \const{X}{s+1} \,;
\]
here and in what follows we use the 
norm on contact diffeomorphisms
$\|\FPX\|_{s+1} := \|\Psi (X)\|_{s+1}$, where 
$\FPX = \exp \circ \Psi (X)$.
\end{remark}

\begin{remark}
\label{rem:standard-estimate}
We will repeatedly use the estimates
\begin{align*}
 \|fg\|_{s} &\prec \|f\|_{s}\|g\|_{s-1} + \|f\|_{s-1}\|g\|_{s} \\
\intertext{and}
 \|g \circ F\|_{s} &\prec \left(\|g\|_{s} + \|g\|_{s}\|F\|_{s-1} +
   \|g\|_{s-1}\|F\|_{s}\right)
	 \cdot \left( 1 + \|F\|_{s-1}\right)^{s-1}\\
 &  \prec \; \|g\|_{s} + \|g\|_{s}\|F\|_{s-1} + \|g\|_{s-1}\|F\|_{s}
\end{align*}
for all $s > 2n+4, f,g \in \Wk{s}(M)$, and $F \in \cdiff{s}, \, \|F\|
\prec 1$, without comment. The first estimate was proved in \cite{BD1}. The
second estimate follows easily by writing $g\circ F$ in local coordinates
and computing $\| g\circ F\|_{U,s}$ in a coordinate neighborhood $U\subset
M$ using the chain rule. In the last estimate, we used the fixed bound on
$X$ to conclude that $(1 + \|F\|_{s-1})^{s-1} \prec 1$.

\end{remark}

Our next goal is to obtain estimates for the deformation tensor for the pull-back $\mu
= \FPX^{\str}\phi$. 
\begin{lemma}
\label{lem:phicircF-estimate}
Let $F = \FPX$ and $s >2n+4$. Then
\[
     \|\phi\circ F\|_{s} \prec 
     \|\phi\|_{s} 
   + \|\phi\|_{s}   \cdot \|X\|_{s-1} 
   + \|\phi\|_{s-1} \cdot \|X\|_{s} \,.
\]
\end{lemma}
\begin{proof}
Observe that on each coordinate patch $U_{\ell}$
\begin{eqnarray*}
\|\rho_{\ell}(\phi  \circ F)_{\ell}\|_{U_{\ell},s}
    &\prec& \|\phi\|_{U_{\ell},s}
	 + \|\phi\|_{U_{\ell},s}\| F_{\ell}\|_{U_{\ell},s-1}
	 + \|\phi\|_{U_{\ell},s-1}\| F_{\ell}\|_{U_{\ell},s} \\
    &\prec & \|\phi\|_{s} + \|\phi\|_{s}\| F\|_{s-1}
	+ \|\phi\|_{s-1}\| F\|_{s}\\
    & \prec& \|\phi\|_{s} + \|\phi\|_{s}\|\Psi( X)\|_{s-1}
	+ \|\phi\|_{s-1}\|\Psi( X)\|_{s}\\
    & \prec& \|\phi\|_{s} + \|\phi\|_{s}\|X\|_{s-1}
	+ \|\phi\|_{s-1}\| X\|_{s}\,
\end{eqnarray*}
The result follows from finiteness of the cover $U_{\ell}$.
\end{proof}

Next let $\cE(X,Y,\phi)$ be the vector-valued one-form defined by
Equation~\eqref{eqn:EX-global}.
Then we have the following estimates:

\begin{lemma}
\label{lem:def-estimate-1}
For $s > 2n+4$, let $\phi\in
\Wk{s}(\Deform)$ be a deformation tensor with $\|\phi\|_{s} <
\const{\phi}{s}$  and let $X, Y \in \Wkcont{s+1}$ be  vector fields with
$\|X\|_{s+1} < \const{X}{s+1}$, $\|Y\|_{s+1} < \const{X}{s+1}$, for $\const{\phi}{s}$ chosen as in 
Remark \eqref{rem:neighbourhood}.
Then
\[
	 \|\cE(X,Y,\phi)\|_{s} \prec 
	 ( \|X\|_{s} + \|\phi\circ \FPY \|_{s} )\cdot \|X\|_{s+1} \,.
\]
Let $\phi_j \in \Wk{s}(\Deform)$, $j=1,2$ be two deformation tensors with
$\|\phi_j\|_{s} < \const{\phi}{s}$, and let $X_{j},Y_{j} \in \Wkcont{s+1}$,
$j=1,2$ be  contact vector fields with $\|X_j\|_{s+1} < \const{X}{s+1}$, $\|Y_j\|_{s+1} < \const{Y}{s+1}$.
Then
\begin{align*}
	 \|\cE(X_1,Y_1,\phi_1) - \cE(X_2,Y_2,\phi_2) \|_{s} \prec& 
  \quad  ( \|X_1\|_{s+1} + \|X_2\|_{s+1} )\cdot \|X_1 - X_2\|_{s}\\
  & + ( \|X_1\|_{s} + \|X_2\|_{s} )\cdot \|X_1 - X_2\|_{s+1}\\
  & + ( \|X_1\|_{s+1} + \|X_2\|_{s+1} )\cdot
	 \|\phi_1\circ F_1 - \phi_2\circ F_2\|_{s}\\
  & + \; (\|\phi_1\circ F_1\|_{s} + \|\phi_2\circ F_2\|_{s}) \cdot
	   \|X_1 - X_2\|_{s+1}
 \,.
\end{align*}
where  $F_j = F_{\Psi(Y_j)}$, $j=1,2$. 
\end{lemma}

\begin{proof}
  By Equation~\eqref{eqn:EX-global}, our proof  amounts to obtaining
  sufficiently good estimates on the entries of the matrices $A(X,Y,\phi)$ and $B(X,Y,\phi)$
  defined in Equations~\eqref{eqn:deformAB}
  Recall the local \formulae\ for $A$ and $B$:
\begin{align*}
A^{\a}_{\b} &=    \delta^{\a}_{\b} + Z_{\b}\inter\Lie_{X}\w^{\a}  +
Z_{\b}\inter 
	 \Quad^{\a}(X,Y,\phi) \\
B^{\a}_{\bb} &= (\bDb X)_{\bb}^{\a} + (\phi^{\a}_{\bb}\circ \FPY) + 
		  Z_{\bb} \inter \Quad^{\a}(X,Y,\phi) \,.
\end{align*}
The estimate 
\begin{align*}
    \|\cE^{\a}_{\bb}\|_{s} 
      &\prec  \|Z_{\bb} \inter \Quad^{\a}(X,Y,\phi)\|_{s} +
      \|A^{-1}\|_{s} \| [(I-A)B]^{\a}_{\bb} \|_{s}
\end{align*}
follows immediately from the formula for $\cE^{\a}_{\bb}$.

We estimate each term on the right-hand side.
First, using Proposition~\ref{prop:Psi-estimates} to estimate
$\Psi(X)-X$ and observing that the estimate
$\norm{Z_{\b}\inter\Quad_{\w^{\a}}(\Psi(X))}_{s} \prec
\norm{\Psi(X)}_{s}\norm{\Psi(X)}_{s+1}$ follows immediately from the local formula
\eqref{eq:Quad-local}, we obtain
\begin{align*}
\label{eqn:cQ-estimate}
\|Z_{\b}\inter \Quad^{\a}(X,Y,\phi) \|_{s} 
\prec&\;
\|Z_{\b}\inter \Lie_{\Psi(X)-X}\w^{\a}\|_{s} +
    \|( \phi^{\a}_{\gb}\circ \FPY)
    Z_{\b}\inter\Lie_{\Psi(X)}\w^{\gb}\|_{s}\\
&  \quad +  
\|Z_{\b}\inter\Quad_{\w^{\a}}(\Psi(X))\|_{s}
 +\| (\phi^{\a}_{\gb}\circ \FPY) Z_{\b}\inter\Quad_{\w^{\gb}}(
 \Psi(X))\|_{s}\\
\prec&\;
\| \Psi(X)-X\|_{s+1} + \|\phi\circ \FPY\|_{s}
\,\|Z_{\b}\inter\Lie_{\Psi(X)}\w^{\gb}\|_{s-1}\\
&\quad +
 \|\phi\circ \FPY\|_{s-1} \,\|Z_{\b}\inter\Lie_{\Psi(X)}\w^{\gb}\|_{s} \\
&\quad +
\|Z_{\b}\inter\Quad_{\w^{\a}}(\Psi(X))\|_{s}
+ \|\phi\circ \FPY\|_{s-1} \|Z_{\b}\inter\Quad_{\w^{\gb}}(\Psi(X))\|_{s}\\
&\quad + \|\phi\circ \FPY\|_{s} \|Z_{\b}\inter\Quad_{\w^{\gb}}( \Psi(X))\|_{s-1}
\\
\prec&\;
\|X\|_{s}\, \|X\|_{s+1} + \|\phi\circ \FPY\|_{s} \,
\|\Psi(X)\|_{s}
+
 \|\phi\circ \FPY\|_{s-1} \,\|\Psi(X)\|_{s+1} \\
&\quad +
  \| \Psi(X)\|_{s}  \, \| \Psi(X)\|_{s+1} 
+ \|\phi\circ \FPY\|_{s-1} \| \Psi(X)\|_{s}\,\| \Psi(X)\|_{s+1}\\
& 
\quad + \|\phi\circ \FPY\|_{s} \| \Psi(X)\|_{s-1} \, \| \Psi(X)\|_{s}  \\
\prec&\;
\|X\|_{s}\, \|X\|_{s+1} + \|\phi\circ \FPY\|_{s} \,
\|X\|_{s}
+
 \|\phi\circ \FPY\|_{s-1} \,\|X\|_{s+1} \\
&\quad +
  \|X\|_{s}  \, \|X\|_{s+1} 
+ \|\phi\circ \FPY\|_{s-1} \|X\|_{s}\,\|X\|_{s+1}
+ \|\phi\circ \FPY\|_{s} \|X\|_{s-1} \, \|X\|_{s}  \\
\prec&\;
 ( \|X\|_{s} + \|\phi\circ \FPY\|_{s} )\cdot \|X\|_{s+1}
\end{align*}
with a similar estimate for $\|Z_{\bb}\inter \Quad^{\a}(X,Y,\phi) \|_{s} $.
Next
\begin{align*}
\| [(I-A)]^{\a}_{\g} \|_{s}
=&\;
\| ( Z_{\g}\inter\Lie_{X}\w^{\a}  + Z_{\g}\inter \Quad^{\a}(X,Y,\phi))  \|_{s}
\\
\prec&\;
\|X\|_{s+1} + 
 ( \|X\|_{s} + \|\phi\circ \FPY\|_{s} )\cdot \|X\|_{s+1} \,,
\end{align*}
which implies in particular that $A = I - (I - A)$ is invertible.  More
precisely, because the matrix $A = \left[A^{\a}_{\b}\right]$ is the of the
form $I + \text{(small matrix)}$, a series expansion for $A^{-1}$ yields
the estimate $\|A^{-1}\|_{s} \prec \|X\|_{s+1} + \|\phi\circ \FPY\|_{s}\cdot
\|X\|_{s+1}$ which is uniformly bounded by a constant depending only on the
constant $C$ in Remark~\ref{rem:neighbourhood}.  Also
\begin{align*}
\|  B^{\g}_{\bb} \|_{s}
=&\;
\| ((\bDb X)_{\bb}^{\g} + (\phi^{\g}_{\bb}\circ \FPY) + 
		  Z_{\bb} \inter \Quad^{\g}(X,Y,\phi)) \|_{s} 
\\
\prec&\;
 \| X\|_{s+1}\ + \|\phi\circ \FPY\|_{s} + 
 ( \|X\|_{s} + \|\phi\circ \FPY\|_{s} )\cdot \|X\|_{s+1} 
 \,,
\end{align*}
so
\begin{align*}
\| [(I-A) B]^{\a}_{\bb} \|_{s}
\prec&\;
\| [(I-A)]^{\a}_{\g} \|_{s} \| [ B]^{\g}_{\bb} \|_{s-1}
+ \| [(I-A)]^{\a}_{\g} \|_{s-1} \| [ B]^{\g}_{\bb} \|_{s}
\\
\prec&\;
\left( \|X\|_{s+1} + 
 ( \|X\|_{s} + \|\phi\circ \FPY\|_{s} )\cdot \|X\|_{s+1} \right) \\
&\qquad \cdot
 \left(  \| X\|_{s}\ + \|\phi\circ \FPY\|_{s-1} + 
 ( \|X\|_{s-1} + \|\phi\circ \FPY\|_{s-1} )\cdot \|X\|_{s} \right)
\,,
 \\
&\quad + \left(\|X\|_{s} + 
 ( \|X\|_{s-1} + \|\phi\circ \FPY\|_{s-1} )\cdot \|X\|_{s} \right)\\
& \qquad \cdot
 \left( \| X\|_{s+1}\ + \|\phi\circ \FPY\|_{s} + 
 ( \|X\|_{s} + \|\phi\circ \FPY\|_{s} )\cdot \|X\|_{s+1} \right)
\,,
\\
\prec&\;
\|X\|_{s+1} \cdot
(\|X\|_{s} +\| (\phi\circ \FPY)\|_{s-1})\\
&\;\quad
+ \|X\|_{s} \cdot
(\|X\|_{s+1} +\| (\phi\circ \FPY)\|_{s})\\
\prec &\;
 ( \|X\|_{s} + \|\phi\circ \FPY\|_{s} )\cdot \|X\|_{s+1} \,.
\end{align*}
This completes the proof of the first estimate.

To prove the second  estimate,
let $A_j = A(X_j,Y_j,\phi_j)$,  $B_j = B(X_j,Y_j,\phi_j)$,
$\cE_j=\cE(X_j,Y_j,\phi_j)$, $j=1,2$.
Then
 \begin{eqnarray*}
	A_{1}^{-1} B_{1} - A_{2}^{-1} B_{2} & = &
	A_{1}^{-1} (B_{1} - B_{2}) - (A_{2}^{-1} - A_{1}^{-1}) B_{2}  \\
	 & = & [(B_{1} - B_{2}) + A_{1}^{-1}  (I - A_{1})(B_{1} - B_{2})] -
         [A_{2}^{-1}
		  (A_{1} - A_{2})A_{1}^{-1} B_{2}] \,.
\end{eqnarray*}
Using this in Equation~\eqref{eqn:EX-global}, we  obtain the equality
\[
\begin{split}
 [\cE_1 - \cE_2]^{\a}_{\bb} &=
  Z_{\bb}\inter \left( \Quad^{\a}(X_1,Y_1,\phi_1) - \Quad^{\a}(X_2,Y_2,\phi_2)
  \right)\\
 &\quad + [ A_{1}^{-1}  (I - A_{1})(B_{1} - B_{2})]^{\a}_{\bb} 
  - [ A_{2}^{-1}   (A_{1} - A_{2})A_{1}^{-1} B_{2} ]^{\a}_{\bb} \,.
\end{split}
\]
Choose the constant $\const{X}{s+2}$ in Remark~\ref{rem:neighbourhood}
sufficiently small to ensure that $\|A_j^{-1}\|_{s}<C'$ for some fixed
constant $C'$.  The triangle inequality, then gives
\begin{multline}
\label{eqn:cauchy-estimate-2}
 \|[\cE_1 - \cE_2]^{\a}_{\bb}\|_{s} \prec
 \| Z_{\bb}\inter \Quad^{\a}(X_1,Y_1,\phi_1)  - Z_{\bb}\inter\Quad^{\a}(X_2,Y_2,\phi_2) \|_{s}\\
  + \|(I - A_{1})(B_{1} - B_{2}) \|_{s} 
  + \| (A_{1} - A_{2})\|_{s} \| B_{2} \|_{s-1} \\
 + \| (A_{1} - A_{2})\|_{s-1} \| B_{2} \|_{s} \,. 
\end{multline}
We  estimate  all four  terms on the right-hand side of
\eqref{eqn:cauchy-estimate-2} in a similar manner. We present the estimate
of the first
term in detail and leave the verification of the estimates of the remaining two terms to the reader.
Rearranging terms and simplifying gives
\begin{align*}
 Z_{\bb}\inter \Quad^{\a}(X_1,Y_1,\phi_1) &- Z_{\bb}\inter \Quad^{\a}(X_2,Y_2,\phi_2)
 \\
=
\quad
&\;
\left\{
   Z_{\bb}\inter\Lie_{\Psi(X_1)-X_1}\w^{\a} +
    ( \phi^{\a}_{1,\gb}\circ F_1)  Z_{\bb}\inter\Lie_{\Psi(X_1)}\w^{\gb} 
\right.\\
&\; 
\left.
   +  
      Z_{\bb}\inter\Quad_{\w^{\a}}(\Psi(X_1)) + (\phi^{\a}_{1,\gb}\circ
      F_1)Z_{\bb}\inter\Quad_{\w^{\gb}} (\Psi(X_1))
\right\}\\
-&\;
\left\{
   Z_{\bb}\inter\Lie_{\Psi(X_2)-X_2}\w^{\a} +
    ( \phi^{\a}_{2,\gb}\circ F_2)  Z_{\bb}\inter\Lie_{\Psi(X_2)}\w^{\gb}
\right. \\
&\;
\left.     
   +  
      Z_{\bb}\inter\Quad_{\w^{\a}}(\Psi(X_2)) + (\phi^{\a}_{2,\gb}\circ
      F_2)Z_{\bb}\inter\Quad_{\w^{\gb}}(\Psi(X_2))
\right\}\\
=
\quad
&\;
   Z_{\bb}\inter\Lie_{(\Psi(X_1)-X_1) - (\Psi(X_2)-X_2)}(\w^{\a} )
+(Z_{\bb}\inter\Quad_{\w^{\a}} (\Psi(X_1)) - Z_{\bb}\inter\Quad_{\w^{\a}}(\Psi(X_2) )
\\
+&\; \left\{
    ( \phi^{\a}_{1,\gb}\circ F_1)  Z_{\bb}\inter\Lie_{\Psi(X_1)}\w^{\gb}
    -
( \phi^{\a}_{2,\gb}\circ F_2)  Z_{\bb}\inter\Lie_{\Psi(X_2)}\w^{\gb} 
\right\}\\
+&\;
\left\{
	(\phi^{\a}_{1,\gb}\circ F_1)Z_{\bb}\inter\Quad_{\w^{\gb}}(\Psi(X_1))   
     - 
      (\phi^{\a}_{2,\gb}\circ F_2)Z_{\bb}\inter\Quad_{\w^{\gb}}( \Psi(X_2) )
\right\}
\end{align*}
By our  previous estimates,  we may estimate as follows:
\begin{align*}
 \|Z_{\bb}\inter \Quad^{\a}(X_1,Y_1,\phi_1) &- Z_{\bb}\inter
 \Quad^{\a}(X_2,Y_2,\phi_2)\|_{s}\\
\prec\quad
&
   \|Z_{\bb}\inter\Lie_{(\Psi(X_1)-X_1) - (\Psi(X_2)-X_2)} \w^{\a}
   \|_{s}\\
&
+\;
\|Z_{\bb}\inter\Quad_{\w^{\a}} (\Psi(X_1)) - Z_{\bb}\inter\Quad_{\w^{\a}}(\Psi(X_2)) \|_{s}
\\
&+\; 
    \|( \phi^{\a}_{1,\gb}\circ F_1)  Z_{\bb}\inter\Lie_{\Psi(X_1)}\w^{\gb}
    -
( \phi^{\a}_{2,\gb}\circ F_2)  Z_{\bb}\inter\Lie_{\Psi(X_2)}\w^{\gb}
\|_{s}
\\
&+\;
\|
	(\phi^{\a}_{1,\gb}\circ F_1)\,Z_{\bb}\inter\Quad_{\w^{\gb}}(\Psi(X_1))
     -  (\phi^{\a}_{2,\gb}\circ F_2)\,Z_{\bb}\inter\Quad_{\w^{\gb}}(\Psi(X_2))
\|_{s}\\
\prec\quad
& \|(\Psi(X_1)-X_1) - (\Psi(X_2)-X_2)\|_{s+1} \\
&+\;
 \| X_1 -X_2\|_{s} \cdot (\| X_1\|_{s+1}+ \|X_2 \|_{s+1}) +\; \| X_1 -X_2\|_{s+1} \cdot (\| X_1\|_{s}+ \|X_2 \|_{s})
\\
&+\; 
    \|( \phi^{\a}_{1,\gb}\circ F_1)  Z_{\bb}\inter\Lie_{\Psi(X_1)}\w^{\gb}
    -
( \phi^{\a}_{2,\gb}\circ F_2)  Z_{\bb}\inter\Lie_{\Psi(X_2)}\w^{\gb}
\|_{s}
\\
&+\;
\|
	(\phi^{\a}_{1,\gb}\circ F_1)Z_{\bb}\inter\Quad_{\w^{\gb}}(\Psi(X_1)   )
     -  (\phi^{\a}_{2,\gb}\circ F_2)Z_{\bb}\inter\Quad_{\w^{\gb}}(\Psi(X_2))
\|_{s} \\
\intertext{where we have used Lemma~\ref{lem:remainder-estimate}(c)}
\prec\quad
&
 \|X_{1} - X_{2}\|_{s}(\|X_{1}\|_{s+1} + \|X_{2}\|_{s+1})
+
 \| X_1 -X_2\|_{s+1} \cdot (\| X_1\|_{s}+ \|X_2 \|_{s}),
\\
&+\; 
    \|( \phi^{\a}_{1,\gb}\circ F_1)  Z_{\bb}\inter\Lie_{\Psi(X_1)}\w^{\gb}
    -
( \phi^{\a}_{2,\gb}\circ F_2)  Z_{\bb}\inter\Lie_{\Psi(X_2)}\w^{\gb}
\|_{s}
\\
&+\;
\|
	(\phi^{\a}_{1,\gb}\circ F_1)Z_{\bb}\inter\Quad_{\w^{\gb}}(\Psi(X_1) )
     -  (\phi^{\a}_{2,\gb}\circ F_2)Z_{\bb}\inter\Quad_{\w^{\gb}}(\Psi(X_2))
\|_{s}, 
\end{align*}
where we have used Proposition~\ref{prop:Psi-estimates}\eqref{cauchy-estimate}.
Observe that
\begin{align*}
    \|( \phi^{\a}_{1,\gb}\circ F_1)& Z_{\bb}\inter\Lie_{\Psi(X_1)}\w^{\gb}
    -
( \phi^{\a}_{2,\gb}\circ F_2)  
Z_{\bb}\inter\Lie_{\Psi(X_2)}\w^{\gb}  \|_{s}\\
\prec\quad&
  (\|\phi_1\circ F_1\|_{s-1}  
+ \|\phi_2\circ F_2\|_{s-1}) \cdot \|X_1 -  X_2\|_{s+1}\\
&\quad +
(\|\phi_1\circ F_1\|_{s} + \|\phi_2\circ F_2\|_{s}) \cdot
\|X_1 - X_2\|_{s} \\
&\quad +\;
(\|X_1\|_{s} + \|X_2\|_{s}) 
   \cdot \|\phi_1 \circ F_1 - \phi_2 \circ F_2\|_{s} \\
&\quad + 
(\|X_1\|_{s+1} + \|X_2\|_{s+1}) 
   \cdot \|\phi_1 \circ F_1 - \phi_2 \circ F_2\|_{s-1}\\
\prec\quad&
(\|\phi_1\circ F_1\|_{s} + \|\phi_2\circ F_2\|_{s}) \cdot \|X_1 -
X_2\|_{s+1}\\
&+\;
(\|X_1\|_{s+1} + \|X_2\|_{s+1}) \cdot \|\phi_1 \circ F_1 - \phi_2 \circ
F_2\|_{s}\,,
\end{align*}
where we have used the identity 
$f_1\, g_1 - f_2\, g_2  = f_1\,( g_1 -  g_2)  + (f_1 - f_2)\, g_2 \;$ and
the corresponding estimate
\begin{align*}
   \| f_1\, g_1 - f_2\, g_2 \|_{s} \prec\quad&
  (\| f_1\|_{s-1} + \| f_2\|_{s-1}) \cdot \| g_1 - g_2\|_{s} +
(\| f_1\|_{s} + \| f_2\|_{s})  \cdot \| g_1 - g_2\|_{s-1}
\\  
+\,&(\| g_1\|_{s-1} + \| g_2\|_{s-1}) \cdot \| f_1 - f_2\|_{s}  +
  (\| g_1\|_{s} + \| g_2\|_{s}) \cdot \| f_1 - f_2\|_{s-1}  \,.
\end{align*}
A similar argument yields the estimate
\begin{multline*}
\|      (\phi_{1}\circ F_1)Z_{\bb}\inter\Quad_{\overline{\w}}(\Psi(X_1))
     -  (\phi_{2}\circ F_2)Z_{\bb}\inter\Quad_{\overline{\w}}(\Psi(X_2))
\|_{s} \\
\prec
\quad
(\|\phi_1\circ F_1\|_{s} + \|\phi_2\circ F_2\|_{s}) 
\cdot \|X_1 - X_2\|_{s+1} \\
+ (\|X_1\|_{s+1} + \|X_2\|_{s+1})
  \cdot \|\phi_1\circ F_1 - \phi_2 \circ F_2\|_{s} \,.
\end{multline*}
Thus
\[
\begin{split}
\|Z_{\bb}\inter \Quad^{\a}(X_1,Y_1,\phi_1) - Z_{\bb}\inter
\Quad^{\a}(X_2,Y_2,\phi_2)\|_{s} 
\prec
  &\quad  ( \|X_1\|_{s+1} + \|X_2\|_{s+1} )\cdot \|X_1 - X_2\|_{s}\\
&  + ( \|X_1\|_{s} + \|X_2\|_{s} )\cdot \|X_1 - X_2\|_{s+1}\\
&  + ( \|X_1\|_{s+1} + \|X_2\|_{s+1} )\cdot
	 \|\phi_1\circ F_1 - \phi_2\circ F_2\|_{s} \\
&  + \; (\|\phi_1\circ F_1\|_{s} + \|\phi_2\circ F_2\|_{s}) \cdot
	   \|X_1 - X_2\|_{s+1} \,.
\end{split}
\]
\nopagebreak[4]\end{proof}

Our proof of the \apriori\ estimates from which
Theorem~\ref{strong-normal-form} follows requires one more technical lemma.
For $k>2n+4$ and $\epsilon>0$ small, let  $\phi \in \Wk{k}(\Deform)$ and $X_0 \in \Wkcont{k}$ with
$\norm{\phi}_{k}<\epsilon$ and $\norm{X_0}<\epsilon$. Then the map
\[
         T_{\phi,X_0}^{k+1} \st \Wkcxcont{k+1} \rightarrow \Wkcxcont{k+1} \st Z \mapsto Z +
         \cPh(\cE(\pi_{\Re}(Z),X_0,\phi) \,,
\]
is defined for all $Z$ in a sufficiently small ball about the origin.

\begin{lemma}
\label{lem:fixed-point} There exists a sufficiently small  $\epsilon>0$
such that the following holds.
For all $\phi\in
\Wk{k}(\Deform)$ and  $X_0 \in \Wkcont{k}$  such that $\norm{X_0}_{k} < \epsilon$ and
$\norm{\phi}_{k}<\epsilon$, the equation $T_{\phi,X_0}^{k+1}(Z) = W$ has a
unique solution $Z \in \Wkcxcont{k+1}$ for all $W\in \Wkcxcont{k+1}$ with
$\norm{W}_{k+1}<\epsilon$. Moreover, the solution satisfies the estimate
$\norm{Z}_{k+1} \leq 2 \, \norm{W}_{k+1}$.
\end{lemma}

\begin{proof}
We first show that we can choose $\delta>0$ so that $Z\mapsto \cPh(\cE(\pi_{\Re}(Z),X_0,\phi)$ is a contraction
mapping in $\Wkcxcont{k+1}$ for $\norm{Z}_{k+1} <\delta$.
To see this, first note that by Lemma~\ref{lem:def-estimate-1}, for $\phi\in
\Wk{k}(\Deform)$ with $\|\phi\|_{k} < C$, for $C$ sufficiently small, 
the estimate 
\[
	 \|\cE(X_1,X_0,\phi) - \cE(X_2,X_0,\phi)\|_{k} \prec 
	 ( \|X_1\|_{k} + \|X_2\|_{k+1} + \|\phi\circ F_0 \|_{k} )\cdot \|X_1-X_2\|_{k+1} 
\]
holds for all $X_1, X_2 \in \Wkcont{k+1}$, with $\|X_j\|_{k+1} <C$, $j=1,2$. Thus,
\[
             \norm{ \cPh(\cE(X_1,X_0,\phi) - \cPh(\cE(X_2,X_0,\phi) )}_{k+1} \prec 	 
            ( \norm{X_1}_{k+1} + \norm{X_2}_{k+1}+ \norm{\phi\circ F_0}_{k}
            )\cdot \norm{X_1 - X_2}_{k+1}  \,.
\]
Consequently, for  $\delta'>0$ sufficiently small, 
\[
             \norm{ \cPh(\cE(X_1,X_0,\phi) - \cPh(\cE(X_2,X_0,\phi)
               )}_{k+1} < \frac{1}{2} \norm{X_1 - X_2}_{k+1}  \,,
\]
provided  $\norm{\phi}_{k}<\delta$, $\norm{X_j}_{k+1} < \delta'$, $j=1,2$.

Now choose $\delta <\delta'$ so that $\norm{\pi_{\Re}(Z)}_{k+1} <\delta'$ for $\norm{Z}_{k+1}<\delta$. 
Then for  $X_j = \pi_{\Re}(Z_j)$, $j=1,2$, 
\[
             \norm{ \cPh(\cE(\pi_{\Re}Z_1,X_0,\phi) - \cPh(\cE(\pi_{\Re}Z_2,X_0,\phi)
               )}_{k+1} < \frac{1}{2} \norm{Z_1 - Z_2}_{k+1}  \,,
\]
provided  $\norm{Z_j}_{k+1} < \delta$.

Finally, set $\epsilon = \delta/2$. Choose any $W \in \Wkcxcont{k+1}$ and define the sequence $Z_n$, $n=0,1,2,\dots$ inductively 
by $Z_0 = 0$, $Z_{n+1} = W - \cPh(\cE(\pi_{\Re}(Z_n),X_0,\phi)$.
Since $\cE(0,X_0,\phi)=0$, $Z_1 = W$. Consequently $\{ Z_n \}$ is Cauchy
with $\norm{Z_{n+1}-Z_{n}}_{k+1} < \frac{1}{2}\norm{Z_n -
    Z_{n-1}}_{k+1}$. Therefore, $\norm{Z_n}_{k+1} < 2 \norm{W}_{k+1}$. Thus, the sequence converges to
a solution $Z$ of the equation $T_{\phi,X_0}(Z) =W$ satisfying the estimate
$\norm{Z}_{k+1} \leq 2 \norm{W}_{k+1}$. Uniqueness of the solution follows from the contraction mapping property.
\end{proof}

We are now able to obtain the \apriori\ estimates that we promised and from
which Theorem~\ref{strong-normal-form} follows.

\begin{theorem}
\label{apriori-estimates-1}
Fix a smooth background CR~structure on $M$ as above, and   let 
\[
\cPh:
 \FormsH{1} \rightarrow
 \Wkcxcont{\infty}
\text{ and } 
{\cQh} : \FormsH{1}
 \rightarrow \FormsH{1}
\]
be the linear operators of Corollary~\ref{cor:homotopy-contact}. 
Then for $s>2n+4$, there exists  $\epsilon > 0$ such that the following holds:

\noindent
Suppose that
$\phi \in \Wk{s+2}(\Deform)$, $X \in \Wkcont{s+1}$ (so $\FPX \in
\cdiff{s+1}$), 
$Y \in V^{s+1}$,
and $\psi \in \Wk{s}(\Deform)$  satisfy the conditions
\[
   \rh (X - i Y ) = 0\,,\quad 
   \|\FPX\|_{s+1} < \epsilon\,,\quad  
   \|\phi\|_{s+2}  < \epsilon\,,
\text{ and }  
    \psi \in \ker\cPh\, .
\]
If the deformation tensor 
$\mu = \FPX^{{\str}}\phi - i \bDb Y - \psi$ is contained in
$\Wk{s+2}(\Deform)$ and $\norm{\mu}_{s+2} < \epsilon$
then 
\[
 \FPX \in\cdiff{s+3}\,,\quad
Y \in V^{s+3}\,\text{ and } \psi \in \Wk{s+2}(\Deform) \,.
\]
Moreover, the following estimates are satisfied:
\begin{align*}
	\|\FPX\|_{s+3}   &\prec  \| \phi\|_{s+2}  + \|\mu\|_{s+2} \,,\\
	\|Y\|_{s+3}      &\prec  \| \phi\|_{s+2}  + \|\mu\|_{s+2}\,, \\
	\|\psi\|_{s+2}   &\prec  \| \phi\|_{s+2}  + \|\mu\|_{s+2}\, .
\end{align*}
\end{theorem}

\begin{proof}  Substitution of the expression for $\FPX^*\phi$ given in
  Proposition~\ref{prop:pull-back-formula} in the formula for $\mu$ gives
\[
          \mu  =  \bDb X
	+  \phi  \circ \FPX  -  i \bDb Y - \psi + \cE(X,X,\phi)
\]
where $\phi\circ F$ and $\cE$ are defined as in
\eqref{eqn:phi-circ-F} and \eqref{eqn:EX-global}.

We first prove   that  $\|X\|_{s+3}$, $\|Y\|_{s+3}$, and
$\|\psi\|_{s+2}$ are finite. 
Applying the operator $\cPh$ and  using the hypothesis 
$\cPh( \psi) = 0$, gives
\begin{equation}
		\cPh(\mu)  =   \cPh(\bDb X -  i \bDb Y)
		+\cPh(\phi  \circ F) + \cPh(\cE(X,X,\phi)) \, .
\end{equation}
Since $\rh (X  - i Y) = 0$, it follows that $X - i Y =
\cPh(\bDb (X -  i Y))$, and solving for
$X - i Y $ in the last equation, we have:
\begin{equation}
\label{eqn:X-solution1}
	 X - i Y  +  \cPh(\cE(X,X,\phi) )=  \cPh(\mu - \phi  \circ \FPX)   \, .
\end{equation}
Next we ``freeze coefficients'' in \eqref{eqn:X-solution1}. Let $X_0-i Y_0
= X-i Y$ and set $W = \cPh(\mu - \phi \circ F_{\Psi(X_0)}) \in
\Wkcxcont{s+3}$.  Then $X_0 - i Y_0$ is the unique solution in
$\Wkcxcont{s+1}$ of the equation
\begin{equation}
\label{eqn:X-solution3}
             T_{\phi,X_0}^{k}(Z) = W\, \text{ for $k=s+1$.}
\end{equation}

We now perform the first of two bootstrapping steps.
Notice that $\phi$ and $\mu$ are small in $\Wk{s+2}$ and, hence, small in
$\Wk{s+1}$, and that  $X_0$ is also small in
$\Wk{s+1}$. Consequently the map $T_{\phi,X_0}^{k}$ is defined for
$k=s+2$. Lemma~\ref{lem:fixed-point} then shows that $T_{\phi,X_0}^{k}$ is
defined for $k=s+2$ and that Equation~\eqref{eqn:X-solution3} with $k=s+2$
has a unique solution in $\Wkcxcont{s+2}$.  It follows that $X_0-i Y_0$ is in
$\Wkcxcont{s+2}$.  Lemma~\ref{lem:phicircF-estimate} then gives the \apriori\ estimate
\[
        \norm{X_0 - i Y_0}_{s+2} \prec \norm{W}_{s+2} \prec \norm{\mu -
          \phi\circ F_{\Psi(X_0)}}_{s+1} \prec \norm{\mu}_{s+1} + \norm{\phi}_{s+1}\,.
\]

The second bootstrap proceeds as follows. We now know that $X_0$ and $\phi$ are both
in $\Wk{s+2}$ and that  $W =  \cPh(\mu - \phi  \circ F_{\Psi(X_0)})$ is in
$\Wk{s+3}$. By shrinking $\epsilon$ if necessary, we can solve
Equation~\eqref{eqn:X-solution3} with $k=s+3$ to conclude that $X_0 - i
Y_0$ is in $\Wkcxcont{s+3}$. Finally, we have $X-i Y = X_0 - i Y_0$ with
the \apriori\ estimate
\[
        \norm{X - i Y}_{s+3} \prec \norm{\mu}_{s+2} + \norm{\phi}_{s+2}\,.
\]
 Finally, since $\psi = F^{{\str}}\phi - i \bDb Y - \mu$, 
it
follows that
$\psi$ is in $\Wk{s+2}$ and satisfies the \apriori\ estimate
\[
\|\psi\|_{s+2} \prec \|\mu\|_{s+2} + \|\phi\|_{s+2} 
+ \|X\|_{s+3} + \|Y\|_{s+3} \prec \|\mu\|_{s+2} + \|\phi\|_{s+2} \,.
\]
This establishes the \apriori\ bounds, and hence the \apriori\ estimates
for $\FPX$, $Y$ and $\psi$.
\end{proof}

\begin{proof}[Proof of Theorem~\ref{strong-normal-form}]
That $F_{\phi}$, $Y_{\phi}$ and $\psi_{\phi}$ are in the appropriate spaces
is an immediate corollary of Theorem~\ref{apriori-estimates-1}.
It remains only to show that the map is continuous.

Choose smooth $\phi_{j}
 \in \Wk{s+2}(\Deform)$ and $\mu_{j} \in \Wk{s+2}(\Deform)$
such that
$\displaystyle
 \phi_{j}\underset{{\scriptscriptstyle \|\cdot\|_{s+2}}}
 {\longrightarrow} \phi
$ and
$\displaystyle \mu_{j}\underset{{\scriptscriptstyle \|\cdot\|_{s+2}}}
 {\longrightarrow} \mu
$.
By the analysis above, there exist
$X_{j} \in \Wkcont{s+3}$,
$Y_{j} \in V^{s+3}$,
and  $\psi_{j} \in \Wk{s+2}(\Harm)$
  such that the contact diffeomorphisms
 $F_{j} = F_{\Psi(X_j)} \in \cdiff{s+3}$ satisfy the
conditions
\[
\mu_{j} = F_{j}^{{\str}}\phi_{j} - i \bDb Y_{j} -   \psi_{j} \,,
\quad
\rh( X_{j} - iY_{j} ) = 0 \,,
\]
with  
$\displaystyle F_{j}\underset{{\scriptscriptstyle \|\cdot\|_{s+1}}}  {\longrightarrow}
F$,
$\displaystyle
Y_{j}\underset{{\scriptscriptstyle \|\cdot\|_{s+1}}} {\longrightarrow} Y
$,
and
$\displaystyle \psi_{j}\underset{{\scriptscriptstyle \|\cdot\|_{s}}} {\longrightarrow}
\psi$.
By the \apriori\ estimates 
$\displaystyle
    \|F_j\|_{s+2} \prec \|\phi_j\|_{s+1} + \|\mu_j\|_{s+1}
$,
$\displaystyle   \|Y_j\|_{s+2} \prec \|\phi_j\|_{s+1} + \|\mu_j\|_{s+1}$,
and
$\displaystyle    \|\psi_j\|_{s+1} \prec \|\phi_j\|_{s+1} + \|\mu_j\|_{s+1}$
established above, we note, in particular, that
$F_j$, $Y_j$, $\psi_j$ are bounded sequences in $\Wk{s+2}$, $\Wk{s+2}$,
and $\Wk{s+1}$, respectively.
Also note that, by continuity of composition,
$F_{j}\underset{{\scriptscriptstyle
    \|\cdot\|_{s+1}}} {\longrightarrow} F$ and
    $\phi_{j}\underset{{\scriptscriptstyle \|\cdot\|_{s+1}}}
    {\longrightarrow} \phi$ together imply 
$\phi_{j}\circ F_{j}\underset{{\scriptscriptstyle
    \|\cdot\|_{s+1}}} {\longrightarrow} \phi\circ F$.
  
We now show that  the sequences $X_j$ and $Y_j$ are Cauchy in
 $\Wkcont{s+2}$.  We estimate $\|X_{j} - X_{i}\|_{s+2}$ as follows.
 Writing
\[
\mu_j  =  \bDb (X_{j} - i \,Y_{j})
	+  \phi_j  \circ F_j  - \psi_j + \cE_j 
\qquad\text{(see \eqref{eqn:phi-circ-F} and \eqref{eqn:EX-global})}
\]
with $F_j = F_{\Psi(X_j)}$, $\cE_j = \cE(X_j,Xj,\phi_j)$
yields the formula
\begin{equation*}
\mu_{j} - \mu_{i}  =  \bDb (X_{j } - X_{i})
	- i \bDb (Y_{j } - Y_{i}) \\
       - (\psi_{j} - \psi_{i}) +
   ( \phi_{j}  \circ F_{j} - \phi_{i}  \circ  F_{i}) + (\cE_{j} -
   \cE_{i})\, .
\end{equation*}
Applying the operator $\cPh$ and  using the facts $\cPh(
\psi_{j}) = 0$,
$\cPh \bDb (X_{j }  - i Y_{j })$ 
$= X_{j }- i Y_{j }$ as above, gives:
\begin{eqnarray*}
\cPh(\mu_{j} - \mu_{i})  &= & \cPh\left(\bDb (X_{j } - X_{i})
  - i \bDb (Y_{j } - Y_{i})\right) 
  + \; \cPh( \phi_{j}\circ F_{j} - \phi_{i}\circ F_{i})
  +  \cPh(\cE_{j} - \cE_{i})
\\
&= &  (X_{j } - X_{i})
  - i  (Y_{j } - Y_{i}) 
  + \;\cPh( \phi_{j}  \circ F_{j} - \phi_{i}  \circ  F_{i})
  +  \cPh(\cE_{j} - \cE_{i}) \,.
\end{eqnarray*}
Solving for $(X_{j } - X_{i}) - i  (Y_{j } - Y_{i}) $, we have:
\[
(X_{j } - X_{i}) - i  (Y_{j } - Y_{i})
     =\cPh(\mu_{j} - \mu_{i})
     - \cPh(\phi_{j}  \circ F_{j} - \phi_{i} \circ F_{i})  - \cPh(\cE_{j} -
     \cE_{i}) \, .
\]
We can estimate the $s+2$ norm for $(X_{j } - X_{i})$ as follows, using
our
\apriori\ bound $\|X_{j}\|_{s+2} \le K$ on the sequence:
\begin{align*}
\|X_{j } &- X_{i}\|_{s+2} + \|Y_{j } - Y_{i}\|_{s+2}\\
&
 \prec
    \|(X_{j } - X_{i}) - i  (Y_{j } -
    Y_{i})\|_{s+2}\\
&
 \prec 
    \|\cPh(\mu_{j} - \mu_{i})\|_{s+2}
  + \|\cPh ( \phi_{j}  \circ F_{j} - \phi_{i} \circ F_{i} )\|_{s+2}  
  + \|\cPh(\cE_{j}- \cE_{i})\|_{s+2}  \\
 &
 \prec
	\|\mu_{j} - \mu_{i}\|_{s+1} + \|\phi_{j}\circ F_{j} - \phi_{i}\circ
        F_{i}\|_{s+1}
    +   \|\cE_{j} - \cE_{i}\|_{s+1}  \\
& \prec 
     \|\mu_{j} - \mu_{i}\|_{s+1} 
   + \|\left( \phi_{j}  \circ F_{j} - \phi_{i}\circ
     F_{i}\right)\|_{s+1}\\
&\qquad
 + (\|X_{i}\|_{s+2}+\|X_{j}\|_{s+2})
	  \cdot \|(\phi_{j}\circ F_{j} - \phi_{i}\circ F_{i})\|_{s+1}
\\
&\qquad
  + (\|\phi_{j}  \circ F_{j}\|_{s+1} + \|\phi_{i} \circ F_{i}\|_{s+1} )
	   \cdot\|X_{j} - X_{i}\|_{s+2} 
\\
&\qquad
   + (\|X_{j}\|_{s+2}+\|X_{i}\|_{s+2})\cdot \|X_{j} - X_{i}\|_{s+1}\\
&\qquad
   + (\|X_{j}\|_{s+1}+\|X_{i}\|_{s+1})\cdot \|X_{j} - X_{i}\|_{s+2}
\end{align*}
\begin{align*}
\phantom{\|X_{j }}&
 \prec
       \|\mu_{j} - \mu_{i}\|_{s+1} 
   +   \|\left( \phi_{j} \circ F_{j} - \phi_{i}\circ F_{i}\right)\|_{s+1}
   +  K \,\|(\phi_{j} \circ F_{j} - \phi_{i} \circ F_{i})\|_{s+1}
\\
&\qquad
   + \left(\|\phi_{j}\|_{s+1} 
       +  \|\phi_{j}\|_{s+1}\|F_{j}\|_{s} + \|\phi_{j}\|_{s}
       \|F_{j}\|_{s+1}
    \right) 
      \cdot \|X_{j} - X_{i}\|_{s+2} 
\\
&\qquad
   + \left(\|\phi_{i}\|_{s+1} 
       +  \|\phi_{i}\|_{s+1}\|F_{i}\|_{s} + \|\phi_{i}\|_{s}
       \|F_{i}\|_{s+1}
    \right) 
      \cdot \|X_{j} - X_{i}\|_{s+2} 
\\
&\qquad
   + K \,  \|X_{j} - X_{i}\|_{s+1}  
   +  \const{}{}\,  \|X_{j} - X_{i}\|_{s+2}
\\
&
 \prec
   \|\mu_{j} - \mu_{i}\|_{s+1} 
	+ \| \phi_{j} \circ F_{j} - \phi_{i}\circ F_{i} \|_{s+1} \\
&\quad +  K \,
   \left(  
       \|(\phi_{j} \circ F_{j} - \phi_{i} \circ  F_{i})\|_{s+1}  
   +   \|    X_{j}             -                  X_{i}\|_{s+1}  
   \right)
 + \const{}{}  \|    X_{j}             -                  X_{i}\|_{s+2}  
\end{align*}
For $ \|\phi\|_{s+1}, \|\mu\|_{s+1}$ sufficiently small (that is
for $\const{}{}$ sufficiently small), we can absorb
the  last  term on the right hand side to obtain an \apriori\ estimate on
the sequence:
\begin{eqnarray*}
 \|X_{j} - X_{i}\|_{s+2} 
& \prec & \|\mu_{j} - \mu_{i}\|_{s+1} 
   +   \|\left( \phi_{j}\circ F_{j} - \phi_{i}\circ F_{i}\right)\|_{s+1}
  +  \|X_{j } - X_{i}\|_{s+1} \, ;
\\
\|\Psi(X_{j}) - \Psi(X_{i})\|_{s+2} &\prec& \|X_{j} - X_{i}\|_{s+2}  \\
&\prec & \|\mu_{j} - \mu_{i}\|_{s+1} 
    +  \|\left( \phi_{j}\circ F_{j} - \phi_{i}\circ F_{i}\right)\|_{s+1} 
  +   \|X_{j } - X_{i}\|_{s+1}\, ;
\\
 \|Y_{j} - Y_{i}\|_{s+2} &\prec & \|\mu_{j} - \mu_{i}\|_{s+1}
       +    \|\left( \phi_{j} \circ F_{j} - \phi_{i}  \circ
         F_{i}\right)\|_{s+1} 
   + \|X_{j } - X_{i}\|_{s+1}\, .
\end{eqnarray*}
Using the facts that $\phi_{j} \circ F_{j}$ , $\mu_j$ ,and $X_j$ are Cauchy
in $\Wk{s+1}$, we have that $X_j$ and $Y_j$ are Cauchy in $\Wk{s+2}$ and
that $X_j\to X$ and $Y_j\to Y$ in $\Wk{s+2}$.

Bootstrapping one more time, using the facts that $\phi_{j} \circ F_{j}$ , $\mu_j$ ,and $X_j$ are Cauchy
in $\Wk{s+2}$, we have that $X_j$ and $Y_j$ are Cauchy in $\Wk{s+3}$ and
that $X_j\to X$ and $Y_j\to Y$ in $\Wk{s+3}$.
This establishes continuity of the map in Theorem~\ref{strong-normal-form}
and completes the proof of the theorem.
\end{proof}



\end{document}